\documentclass[final]{article}

\usepackage{amsfonts}
\usepackage{amssymb}
\usepackage{textcomp}
\usepackage{latexsym,amsmath}
\usepackage{graphics}

\usepackage{hyperref}

\renewcommand{\phi}{\varphi}
\newcommand{\be}{\begin{equation}}
\newcommand{\ee}{\end{equation}}
\newcommand{\ba}{\begin{eqnarray}}
\newcommand{\ea}{\end{eqnarray}}
\newcommand{\ban}{\begin{eqnarray*}}
\newcommand{\ean}{\end{eqnarray*}}

\title{Symmetric interpolatory dual wavelet frames
\thanks{This research was supported by the RFBR-grant (\#12-01-00216-a) and by the grant of the St.Petersburg State University (\#9.38.62.2012).}}

\author{
A.V.~Krivoshein
}
\date{Department of Applied Mathematics and Control Processes,
 \\ St. Petersburg State University \\
 KrivosheinAV@gmail.com, a.krivoshein@spbu.ru}

\begin{document}
\maketitle


\newtheorem{theo}{Theorem}
\newtheorem{lem}[theo]{Lemma}
\newtheorem {prop} [theo] {Proposition}
\newtheorem {coro} [theo] {Corollary}
\newtheorem {defi} [theo] {Definition}
\newtheorem {rem} [theo] {Remark}
\newtheorem {ex} [theo] {Example}

\begin{abstract}
For any symmetry group $H$ and any appropriate matrix dilation we give an explicit method for the construction of $H$-symmetric refinable interpolatory refinable masks  which satisfy sum rule of arbitrary order $n$. For each such mask we give an explicit technique for the construction of dual wavelet 
frames such that the corresponding wavelet masks are mutually symmetric and have the vanishing moments up to the order n. For an abelian symmetry group $H$ we modify the technique such that each constructed wavelet mask is $H$-symmetric.
\end{abstract}




\section{Introduction}

The symmetry of refinable and wavelet functions is one of the
important feature in applications. A lot of papers are devoted to the construction of wavelets with the symmetry properties. In the univariate case the general methods providing the symmetry for wavelets  were presented in~\cite{Han09},~\cite{Pet} (in~\cite{HanZh},~\cite{Zhuang} for multiwavelets) and the references therein.
But in the multivariate case to the best knowledge of the author there are no general methods for the construction of wavelet systems with an arbitrary symmetry properties and for any appropriate dilation matrices. 
The problem in this case is complicated by the fact that there are different kinds of symmetry. 
Also, unlike the univariate case,
in the multivariate case there are no simple and explicit general algorithms for deriving the wavelet systems from the refinable functions. The variety of methods for the construction of wavelets in the different setups in the multivariate case can be found in~\cite{EhlerHan},~\cite{HanShen},~\cite{KrSk}, etc. (and the reference therein).

It seems impossible to provide a complete list of references about the construction of symmetric wavelets.  
The following synopsis just briefly highlight some investigations on the problem. 
The symmetry properties of functions can be defined using the notion of symmetry group (for details see Section 2.1 and, also,~\cite{Han3}). Let $H$ be a symmetry group.
In~\cite{And} G. Andaloro  and co-authors suggest  a method 
for the construction of $H$-symmetric wavelet systems starting with refinable interpolatory masks. The method is based on the lifting scheme; the choice of dilation matrix $M$ is restricted only to $M=2I_d.$
K.~Koch in~\cite{Koch} give an algorithm for the construction of  interpolating scaling vectors with the $H$-symmetry property. 
 The multiwavelets are also constructed using certain technique that does not care about the symmetry properties. In several examples it turns out that the multiwavelets are mutually symmetric. In this paper we clarify this effect in the case of scalar wavelets.
     Some schemes for the construction of highly symmetric wavelet systems were presented 
     by Q. T. Jiang and D.K. Pounds in~\cite{Jiang6fold}.
In~\cite{KST} M. Skopina  and co-authors suggest a method for the construction of refinable interpolatory masks and wavelet masks which are point-symmetric (point-antisymmetric)  for several dilation matrices. In~\cite{Kr} the point and axial symmetric  frame-like wavelet systems were constructed.
The aim of this paper is to extend the results in~\cite{KST} and to generalize the method in~\cite{Kr} on any symmetry group and any appropriate dilation matrix in order to construct the $H$-symmetric dual wavelet frames generated by interpolatory refinable masks. 

 Besides the symmetry property, there are another important features for applications such as the good approximation properties of wavelet system, the smoothness and the small support of refinable and wavelet functions. 
 The good approximation property of wavelet system is connected with the orders of vanishing moments of the wavelet functions, which in turn are connected with the orders of sum rule of the refinable masks. The desirable order of vanishing moments for the wavelet functions will be explicitly incorporated in the proposed construction. The order of smoothness of the refinable functions does not incorporated in the construction. 
 But the Sobolev smoothness exponent can be computed ex post, using the efficient algorithm developed by B. Han in~\cite{HanSymSmooth}. The  analysis of smoothness and its dependence from the mask support  can be found in~\cite{HanAnalSmooth}. Some facts about how the symmetry affects on the smoothness can be found in~\cite{HanAnalSmooth2}.

The paper is organized as follows.  In Section 2 we give some notations and the notion of symmetry group.
In Section 3 we recall the Mixed Extension Principle. 
The main results are contained in Section 4. Subsection 4.1 is devoted to  the construction of interpolatory $H$-symmetric refinable masks. 
    Subsection 4.2 describes the construction of
   mutually symmetric wavelets.
    In Subsection 4.3 we present a method for the construction of wavelet masks such that all of them have the $H$-symmetry properties.
    In Section 5 several examples are presented.


\section{Notation}

${\mathbb N}$ is the set of positive integers,
    ${\mathbb R}$  is the set of real numbers,
    ${\mathbb C}$ is the set of complex numbers.
    ${\mathbb R}^d$ is the
    $d$-dimensional Euclidean space,  $x = (x_1,\dots, x_d)$ and $y =
    (y_1,\dots, y_d)$ are its elements (vectors),
    $(x)_j=x_j$ for $j~=~1,\dots,d,$
    $(x, y)~=~x_1y_1+~\dots~+x_dy_d$,
    $|x| = \sqrt {(x, x)}$, ${\bf0}=(0,\dots, 0)\in{\mathbb R}^d$.
    For $x, y \in{\mathbb R}^d$ we write
    $x > y$ if $x_j > y_j,$ $j= 1,\dots, d;$
    $x \geq y$ if $x_j \geq y_j,$ $j= 1,\dots, d.$
    ${\mathbb Z}^{d}$ is the integer lattice
    in ${\mathbb R}^d$, ${\mathbb Z}_+^d:=\{x\in{\mathbb Z}^{d}:~x\geq~{\bf0}\}.$
    If $\alpha,\beta\in{\mathbb Z}^{d}_+$, $a,b\in{\mathbb R}^d$, we set
    $[\alpha]=\sum\limits_{j=1}^d \alpha_j$,
    $\alpha!=\prod\limits_{j=1}^d\alpha_j!$,
    $\binom{\beta}{\alpha}=\frac{\alpha!}{\beta!(\alpha-\beta)!}$,
    $a^b=\prod\limits_{j=1}^d a_j^{b_j}$,
    $D^{\alpha}f=\frac{\partial^{[\alpha]} f}{\partial x^{\alpha}}=\frac{\partial^{[\alpha]} f}{\partial^{\alpha_1}x_1\dots
    \partial^{\alpha_d}x_d}$, $\delta_{ab}$~is the Kronecker delta.
    Suppose $n\in{\mathbb N}$, then $\Delta_n~:=~\{\alpha:~\alpha\in{\mathbb Z}^{d}_+,\,\,[\alpha]<n\}.$
    Suppose $\beta\in{\mathbb Z}^{d}_+$, then $\square_{\beta}:=\{\alpha: \alpha\in{\mathbb Z}^{d}_+,\,\alpha\le\beta\}.$
    By $\# H$ denote the cardinality of the set $H.$
        By $\texttt{diag}(u_1,\dots,u_d),$ $u_1,\dots,u_d\in{\mathbb Z},$ denote the $d\times d$ diagonal matrix
    such that numbers $u_1,\dots,u_d$ are placed on the main diagonal.
The (i,j)-th element of matrix $M$ is denoted by $[M]_{i,j}$

An  integer  $d\times d$ matrix $M$ is called a  dilation matrix if
the eigenvalues of $M$ are bigger than 1 in module.
    By $M^*$ denote the transpose
 conjugate matrix $M$, by $I_d$ denote the $d\times d$ identity matrix, $m=|\det M|.$ We say
    that the vectors $\alpha,\,\beta\in{\mathbb Z}^{d}$ are congruent modulo $M$ if $\alpha-\beta=M\gamma,$ $\gamma\in{\mathbb Z}^{d}$.
    The integer lattice ${\mathbb Z}^{d}$ is
    split into cosets with respect to the congruence. It is known that the number of the cosets
    is equal to $m$ (see, e.g., \cite[\S~2.2]{NPS}).
    Let us choose
    an arbitrary representative from each coset, call them all digits, and
    denote by $D(M)=\{s_0,\dots, s_{m-1}\}$ the set of digits.
    Thus, any $\alpha\in{\mathbb Z}^{d}$ can be uniquely represented as $\alpha=M\beta+s_k,$
    where $\beta\in {\mathbb Z}^{d},$ $s_k\in D(M).$ 
    Throughout the paper we assume that $s_0={\bf0}.$
    Denote by $\langle s\rangle$ the coset containing digit $s\in D(M)$ or, equivalently, $\langle s\rangle:=\{M\beta+s, \beta\in{\mathbb Z}^{d}\}.$


Let $f$ be a function defined on ${\mathbb R}^d$; then
        $
        f_{jk}=m^{j/2}f(M^j\cdot+k),$ where $j\in{\mathbb Z}, k\in{\mathbb Z}^{d},
        $
    and $\{ f_{jk}\}:=\{f_{jk}:  j\in{\mathbb Z},\,  k\in{\mathbb Z}^{d}\}.$
   The Fourier transform of $f\in L_1({\mathbb R}^d)$
    is defined to be $\widehat
    f(\xi)=\int\limits_{{\mathbb R}^d} f(x)e^{-2\pi i (x, \xi)}\,dx,$ $\xi\in{\mathbb R}^d.$
    This notion can be naturally extended to $L_2({\mathbb R}^d)$.

Let $S'$ denote
    the space of tempered distributions.
    For any $f\in S',$ its
    \emph{Sobolev smoothness exponent} is defined to be

        $$\nu_2(f)=\sup \left\{ \nu\in{\mathbb R}\bigcup\{-\infty\}\bigcup\{+\infty\}: \int_{{\mathbb R}^d} |\widehat f (\xi)|^2 (1+ |\xi|^2)^{\nu} d \xi < \infty\right\}.$$

A trigonometric polynomial $t(\xi),$ $\xi\in{\mathbb R}^d,$ is a finite linear combination of complex exponentials, i.e.  $t(\xi)=\sum_{k\in{\mathbb Z}^{d}}h_k e^{2\pi i \left(k,\xi\right)},$ $h_k\in{\mathbb C}.$  
A shifted trigonometric polynomial $t'(\xi),$ $\xi\in{\mathbb R}^d,$ is defined by $t'(\xi)=e^{2\pi i (\sigma,\xi)} t(\xi)$ $\sigma\in{\mathbb R}.$

A function/distribution $\phi$ is called \emph{refinable} if there exists a $1$-periodic
    function $m_0\in L_2([0,1]^d)$ (refinable mask, also low-pass filter) such that
        \be
        \widehat\phi(\xi)=m_0(M^{*-1}\xi) \widehat\phi(M^{*-1}\xi).\label{RE}
        \ee
    This condition is called the refinement equation.
    It is well
    known (see, e.g.,~\cite[\S~2.4]{NPS}) that for any  trigonometric polynomial
    $m_0$ satisfying $m_0({\bf0})=1$ there exists a unique (up to a factor) solution  of
    the refinement equation~(\ref{RE}) in $S'$. The solution is compactly supported
    and it is given by its Fourier transform
        \be
        \widehat\phi(\xi):=\prod_{j=1}^\infty m_0(M^{*-j}\xi).
        \label{Prod}
        \ee
    Throughout the paper we assume that any refinable mask $m_0$ is a trigonometric
    polynomial with real or complex coefficients and $m_0({\bf0})=1$.

Let us fix a dilation  matrix $M$ and its digits. For any trigonometric polynomial $t$ there exists a unique set of
    trigonometric polynomials $\tau_k$, $k=0,\dots, m-1$, such that
        \be
        t(\xi)= \frac1{\sqrt m}\sum\limits_{k=0}^{m-1} e^{2\pi i(s_k,\xi)}\tau_k(M^*\xi),
        \label{PR}
        \ee
    where $s_k$ are the digits of $M$.
    Equality~(\ref{PR}) is the \emph{polyphase representation} of  $t$.
    For each $k=0,\dots,m-1$ the trigonometric polynomial $\tau_k$ is called the \emph{polyphase component} of $t$
    corresponding to the digit $s_k.$
    Changing the set of digits  will also change
    the polyphase components of the trigonometric polynomial $t$ (namely, they will be shifted by certain integer shifts).
        Let us remark that the polyphase components can be explicitly defined by 
(see, e.g.,~\cite[\S~2.6]{NPS})
        \be
    \tau_{ k}(\xi)=\frac1{\sqrt m}
    \sum\limits_{s\in D(M^{*})} e^{-2\pi i(M^{-1}s_k,\xi+s)}t(M^{*-1}(\xi+s)),
    \quad k=0,\dots, m-1.
        \label{fPolyReprInv}
        \ee

Let  $t$ be a trigonometric polynomial. We say that $t$ is \textit{interpolatory} if 
$$\sum_{s\in D(M^*)}t(\xi+M^{*-1}s)\equiv 1.$$     
It follows from~(\ref{fPolyReprInv}) that for the interpolatory  trigonometric polynomial $t$ its polyphase component corresponding to the digit $s_0={\bf0}$ is a constant, namely $\tau_{0}\equiv\frac 1 {\sqrt{m}}.$

Let $n\in{\mathbb N}.$
We say that a trigonometric polynomial $t$ has \emph{sum rule of order $n$}
    with respect to the dilation matrix $M$ if
        \be
        D^\beta t({M^*}^{-1}\xi)\Big|_{\xi=s}=0, \quad \forall s\in D(M^*)\setminus\{{\bf0}\},
        \quad \forall  \beta\in\Delta_n.
        \label{fSR}
        \ee
    The order of sum rule is a very
    important feature in applications, since it is connected with the good approximation properties of the
    corresponding wavelet systems (see, e.g.~\cite{Han11lpm} and the references therein). Condition~(\ref{fSR}) can be formulated via the polyphase components of  trigonometric polynomial $t$. 
        \begin{lem}~(\cite[Theorem 11]{DSk})
    	A trigonometric polynomial $t$ has sum rule of order $n\in{\mathbb N}$ if and only if its polyphase components $\tau_{k},$ $k=0,\dots, m-1,$ satisfy
        \ba
        D^\beta\tau_{k}({\bf0})=\frac {(2\pi i)^{[\beta]}} {\sqrt m}
        \sum\limits_{\alpha\in\square_\beta}\lambda_\alpha
        \binom{\beta}{\alpha}(- M^{-1}s_k)^{\beta-\alpha} \quad  \forall
        \beta\in\Delta_n,
        \label{fPolySumRule}
        \ea
       where $\lambda_{\alpha}\in \mathbb C$ are defined by
       \ba
       \lambda_\alpha = \frac 1 {(2\pi i)^{[\alpha]}}
        D^{\alpha} t(M^{*-1}\xi)\Big|_{\xi={\bf0}}, \quad
        \alpha\in\Delta_n.
		\label{fPolyphLam}        
        \ea
        \label{lemDynSk}
    \end{lem}
     Note that if trigonometric polynomial $t$ has sum rule of order 1 and $t({\bf 0})=1$ then condition~(\ref{fPolySumRule}) is just
 $\tau_{k}({\bf 0})=\frac 1 {\sqrt m},$ $k=0,\dots, m-1,$ $\lambda_{\bf 0}=1.$
    
For an interpolatory  trigonometric polynomial $t$ Lemma~\ref{lemDynSk} can be simplified.
Since $\tau_{0}=\frac 1 {\sqrt{m}}$, then from condition~(\ref{fPolySumRule}) 
for $\tau_{0}$ we get that $\lambda_\alpha=\delta_{{\bf0}\alpha},$ $\alpha\in\Delta_n.$ 
 Thus, the conditions on the polyphase components of $t$ are
     \ba
     D^{\beta}\tau_{k}({\bf0})=
     \frac {(2\pi i)^{[\beta]}}{\sqrt m}(-M^{-1}s_k)^{\beta}
     \quad \forall
        \beta\in\Delta_n,
    \label{fPolySRInterp}
        \ea 
        where $k=0,\dots,m-1.$  
Also by~(\ref{fPolyphLam}) we have that $D^{\alpha} t(M^{*-1}\xi)\Big|_{\xi={\bf0}}=\delta_{{\bf0}\alpha}$, $\alpha\in\Delta_n.$
    
 We say that a trigonometric polynomial $t$ has \emph{vanishing moments of order $n$} if
        $$D^{\beta} t(\xi)\Big|_{\xi={\bf 0}}= 0, \quad \forall  \beta\in\Delta_n.$$ 
Finally, we say that $t$
has \emph{linear-phase moments of order $n$} with  phase $c\in{\mathbb R}^d$ if
        $$D^{\beta} t({\bf0})=
         D^{\beta} e^{2\pi i (c,\xi)}\Big|_{\xi={\bf0}}= (2\pi i)^{[\beta]} c^\beta,
         \quad \forall  \beta\in\Delta_n.$$
    This notion was described in~\cite{Han10ComplSym} for the univariate case and in~\cite{Han11lpm}
    for more general settings. The importance of the  linear-phase moments is
    in the fact that they are useful in the setting
    of polynomial reproduction and subdivision schemes.

Notice that the above properties of trigonometric polynomials
    are invariant with respect to the integer shifts. If  a trigonometric polynomial $t$
     has sum rule of order $n$ (or vanishing moments of order $n,$ or linear-phase moments of order $n$ with phase $c\in{\mathbb R}^d$)
    then $t'(\xi):=t(\xi) e^{2\pi i (\gamma,\xi)},$ $\gamma\in{\mathbb Z}^{d},$ also has sum rule of order $n$
    (or vanishing moments of order $n,$ or linear-phase moments of order $n$ but with phase $c+\gamma$).

For functions $\psi^{(\nu)}$, $\nu=1,\dots, r$, the system of their dilations and translations
    $\{\psi_{jk}^{(\nu)}\}$ is called a wavelet system.
    Suppose $n\in{\mathbb N}.$ We say that the wavelet system
    $\{\psi_{jk}^{(\nu)}\}$ has
    \emph{vanishing moments of order $n$}
    (or has a {\em  $VM^{n}$ property})
    if
        $$D^{\beta}\widehat{\psi^{(\nu)}}({\bf0})=0, \quad \nu=1,\dots, r, \quad \forall \beta\in\Delta_n.$$

A family of functions $\{f_\alpha\}_{\alpha\in\aleph}$ ($\aleph$ is a countable index set) in a
    Hilbert space $H$ is called a frame in $H$
    if there exist constants $A, B > 0$ such that  
        $$A\|f\|^2\le\sum\limits_{\alpha\in\aleph}|\langle f, f_\alpha\rangle|^2\le B \|f\|^2,\quad \forall f\in H. $$
    The important property of any frame $\{f_\alpha\}_\alpha$ in $H$
    is the following: every $f\in H$ can be decomposed as
    $ f=\sum_\alpha\langle f,\widetilde f_\alpha\rangle f_\alpha, $ where $\{\widetilde f_\alpha\}_\alpha$ is a
    dual frame in $H$. The comprehensive characterization of  frames
    can be found in~\cite{Chris}. The wavelet frames are of great interest in many applications, especially
    in signal processing.

Further, we give some results from the theory of permutation groups (see, e.g.,~\cite{Cameron}). Let $H$ be a finite group (with a binary operation "$\cdot$") and let $\Omega$ be a finite set. A group  action of $H$ on $\Omega$ is a map $\chi: H\times\Omega \longrightarrow \Omega$  such that the following conditions hold for any element $\omega\in\Omega:$ 
    \begin{enumerate}
      \item the identity, i.e. $\chi(I,\omega)=\omega$, where $I$ is the identity element of $H$;
      \item the associativity, i.e. $\chi(E,\chi(\widetilde E,\omega))=\chi(E \cdot\widetilde E,\omega)$ for all $E, \widetilde E \in H$. 
      \end{enumerate}    
The set $\Omega$ is called an $H$-space if there exist a group action $\chi$ of $H$ on $\Omega$.  
A group action of $H$ on $\Omega$ somehow permutes the elements of $\Omega$. 

Suppose $\Omega$ is an $H$-space and $\omega\in\Omega.$ 
 Denote the map $\chi(E,\omega)$ by $E\omega.$
The set $H\omega:=\{E\omega, E\in H\}$ is called the orbit of  $\omega$. Not hard to see that two orbits are either equal or disjoint. Thus, $\Omega$ can be represented as the union of disjoint orbits.

\begin{lem}
Suppose $\Omega$ is an $H$-space. Then there exists a set $\Lambda\subset \Omega$ such that  $\Omega=\cup_{\omega\in \Lambda}H\omega$ and the orbits $H\omega, $ $\omega\in \Lambda,$ are mutually disjoint.
    \label{lemCameron1}
    \end{lem}
 
 The subgroup $H_{\omega}$ of $H$ defined by $H_{\omega}:=\{E\in H: E\omega=\omega\}$ is called the stabilizer of $\omega$. $H\backslash H_{\omega}$ is the quotient group of $H$ modulo $H_{\omega}$. The orbit of $\omega$ is  isomorphic to  $H\backslash H_{\omega}$. Denote by $\Gamma_{\omega}$ a complete set of representatives of the cosets of $H\backslash H_{\omega}$.
 Then, the orbit of $\omega$ can be represented as $H\omega=\cup_{E\in\Gamma_{\omega}}{E\omega}.$ 
 Hence, the set $\Omega$ can be decomposed as 
 \be
	 \Omega=\cup_{\omega\in \Lambda} \cup_{E\in\Gamma_{\omega}}{E\omega},	
 \label{fOmega}
 \ee
 where the set $\Lambda$ is from Lemma~\ref{lemCameron1}.
 Also note that $\#\Gamma_{\omega}=\#H\omega$, the set $\Gamma_{\omega}$ and $H_{\omega}$ consists of the elements of $H$, the set $H\omega$ consists of the elements of  $\Omega.$

    
     \begin{lem} Suppose $\Omega$ is an $H$-space, $\omega\in\Omega.$
Let $\Gamma_{\omega}$ be a complete set of representatives of the cosets of $H\backslash H_{\omega}$. Then
$\#H=\#H_{\omega}\cdot\#\Gamma_{\omega}$
and each $\widetilde E\in H$ can be uniquely represented as $\widetilde E=E\cdot F,$ where $E\in \Gamma_{\omega}$, $F\in H_{\omega}.$
    \label{lemCameron2}
    \end{lem}
 
 Suppose $\omega\in\Omega,$ $\eta\in H\omega$, i.e. there exist $E\in H$ such that $\eta=E\omega.$ Then the stabilizer of $\eta$ and the stabilizer of $\omega$ are conjugate subgroups, namely $E H_{\omega} E^{-1}=H_{\eta}.$

\subsection{Symmetry property of masks and refinable functions}

A finite set $H$ of $d\times d$ unimodular matrices (i.e. square integer matrices with determinant $\pm 1$) is a
    symmetry group on ${\mathbb Z}^{d}$ if $H$ forms a group under the matrix multiplication.
Let $H$ be a symmetry group on ${\mathbb Z}^{d}$.
    A function $f$ is called \emph{$H$-symmetric} with respect to the
    center $C\in {\mathbb R}^d$ if
        $$f(E(\cdot-C)+C)=f, \quad\forall E\in H, \quad x\in{\mathbb R}^d.$$
    For trigonometric polynomials we use a bit another definition which is compatible with the above definition of $H$-symmetric functions. We say that
        $t(\xi)=\sum_{k\in{\mathbb Z}^{d}}h_k e^{2\pi i (k, \xi)},$ $h_k \in{\mathbb C},$
    is an \emph{$H$-symmetric }
    with respect to the center $c\in {\mathbb R}^d$ trigonometric polynomial if
        \be
        t(\xi)=e^{2\pi i (c-Ec,\xi)}t(E^*\xi), \quad
        \forall E\in H
        \label{TrigSym}
        \ee
    and $c-Ec\in{\mathbb Z}^{d}, \forall E\in H.$
    Condition~(\ref{TrigSym}) is equivalent to
        $h_k=h_{E(k-c)+c}, 
        $ $\forall k\in{\mathbb Z}^{d},$ $\forall E\in H.$

Firstly, we show that an $H$-symmetric interpolatory trigonometric polynomial $t$ has some restrictions on the symmetry center. 


\begin{lem}
Let $t$ be an interpolatory trigonometric polynomial that is $H$-symmetric with respect to the center $c\in{\mathbb R}^d$ and has sum rule of order $n$, $n\ge 2$.
Then  $c=Ec$ for all $E\in H.$
\label{lemInterp}
\end{lem}

{\bf Proof.} 
From~(\ref{fPolyphLam}) with $\lambda_\alpha=\delta_{{\bf0}\alpha},$ $\alpha\in\Delta_n,$  and the higher chain formula for the linear change of variables it follows that
\be
        D^{\beta} t(\xi)\Big|_{\xi={\bf0}}=D^{\beta} t(E^*\xi)\Big|_{\xi={\bf0}}=(2\pi i)^{[\beta]}\delta_{{\bf0}\beta}\quad \forall \beta\in\Delta_n.
\label{fMaskIntepr}
\ee
Since $t$ is $H$-symmetric, we have that $\forall E\in H$
and $\forall \beta\in\Delta_n$
        $$D^{\beta} t(M^{*-1}\xi)\Big|_{\xi={\bf0}}=
        D^{\beta} e^{2\pi i (c-Ec,M^{*-1}\xi)}t(E^*M^{*-1}\xi)\Big|_{\xi={\bf0}}=
       (2\pi i)^{[\beta]}
        (M^{-1}(c-Ec))^\beta.
        $$      
Thus, $M^{-1}(c-Ec)={\bf0},$ or, equivalently, $c=Ec,$ $\forall E\in H.$
$\Diamond$

The condition  $c=Ec$ for all $E\in H$ means that either $c={\bf0}$ or $\det(I_d-E)=0$ for all $E\in H$. In the latter case $c$ belongs to the intersection of the null spaces of matrices $I_d-E$, $\,\forall E\in H$ and condition~(\ref{TrigSym}) reduces to $ t(\xi)=t(E^*\xi),$ $\,\forall E\in H.$ Thus, we can say that $t$ is $H$-symmetric with respect to any point from the intersection of the null spaces of matrices $I_d-E,$ $\,\forall E\in H$. So, it is no matter how the symmetry center will be chosen. For example, if      $H:=\left\{ I_2, \left(\begin{matrix}
            -1 & 0\cr
            0 & 1\cr
        \end{matrix}\right) \right\},$
        then the symmetry center of 
        trigonometric polynomial defined in Lemma~\ref{lemInterp}
        can be at any point on the $y-$axis.
For convenience, we assume that the symmetry center for $H$-symmetric interpolatory  trigonometric polynomials is $c={\bf0}.$
Note that the trigonometric polynomial from Lemma~\ref{lemInterp} automatically has linear-phase moments (with phase $c={\bf0}$) equal to the order of sum rule. 

It is known (see, e.g.~\cite{HanSym02}) that the $H$-symmetry of  refinable mask not always carries over
    to its refinable function. Due to this fact
    the notion of symmetry group was modified in~\cite{HanSym02}.

A finite set $H_M$ of $d\times d$ integer matrices  is called
    a \emph{symmetry group with respect to the dilation matrix $M$} if
    $H_M$ is a symmetry group on ${\mathbb Z}^{d}$ such that
$M^{-1}EM\in H_M,$ $\forall E\in H_M.$
   The last property can be also interpreted as follows: for each $E\in H_M$ there exist
   $E'\in H_M$ such that
    \be
    EM=ME' \quad \texttt{or} \quad M^{-1}E=E'M^{-1}.
    \label{fMEEM}
    \ee

The following statement was shown by B.~Han.

\begin{lem}~\cite[Proposition 2.1]{Han3}
    Let $H_M$ be a symmetry group with respect to
    the dilation matrix $M$ and let  $m_0$ be a refinable mask. Then the mask $m_0$ is $H_M$-symmetric with respect to
   the center $c\in{\mathbb R}^d$ if and only if the corresponding refinable function $\phi$ given by~(\ref{Prod})
    is $H_M$-symmetric with respect to the center $C\in {\mathbb R}^d$, where
    $C=(M-I_d)^{-1}c.$
    \label{HanSym}
    \end{lem}

 Let $M$ be a dilation matrix, $H$ be a symmetry group with respect to
    the dilation matrix $M$, $c\in{\mathbb R}^d$ be an appropriate symmetry center, i.e.  $c-Ec\in{\mathbb Z}^{d},$ $\forall E\in H.$ 
It is known that each vector $\beta\in{\mathbb Z}^{d}$ can be uniquely represented as follows: $\beta=M\gamma + s,$ where $\gamma\in{\mathbb Z}^{d}$ and $s\in D(M).$ 
Using this fact it is easy to see that for each digit $s\in D(M)$
 and matrix $E\in H$ there exist a unique digit $q\in D(M)$ and a vector $r_s^E\in {\mathbb Z}^{d}$ such that 
 \be
 Es=M r_{s}^E+q+Ec-c.
 \label{fDigMain}
	\ee 
The indices of $r_s^E$ mean that the vector $r_s^E$ depends on digit $s$ and matrix $E.$

Denote by ${\cal D}:=\{\langle s_0\rangle,...,\langle s_{m-1}\rangle\}$ the set of cosets, where the coset is defined by $\langle s_i\rangle:=\{M\beta+s_i, \beta\in{\mathbb Z}^{d}\},$ $i=0,...,m-1.$   
   Define the map $\chi$ from the set $H\times{\cal D}$  
 as follows 
 	\be
 	\chi(E,\langle s\rangle)=E\langle s\rangle:=\{EM\beta+Es+c-Ec, \beta\in{\mathbb Z}^{d}\}, \quad E\in H.
	\label{fEs} 	
 	\ee
 	
	\begin{prop}
		The set ${\cal D}$ is an H-space, where the group action is defined by~(\ref{fEs}).
    \label{propHspace}
    \end{prop} 	
    
\textbf{    Proof.} Let us fix  $E\in H$ and show that $E\langle s\rangle\in{\cal D}$ for all $s\in D(M).$ 
   Due to~(\ref{fDigMain}) there always exists a unique digit $q\in D(M)$ such that $Es+c-Ec=M\gamma+q,$  $\gamma \in {\mathbb Z}^{d}.$
    Therefore 
    $$E\langle s\rangle=\{EM\beta+Es+c-Ec, \beta\in{\mathbb Z}^{d}\}=
    \{M(E'\beta+\gamma)+q, \beta\in{\mathbb Z}^{d}\}=\langle q\rangle,$$ 
    where  $E'$ is from~(\ref{fMEEM}).  The identity condition in the definition of group action is obviously valid. It remains to show the associativity.
   Suppose $E, \widetilde E\in H.$ Then 
   $$\widetilde E E \langle s\rangle=
   \{\widetilde E EM\beta+\widetilde E Es+c-\widetilde E Ec, \beta\in{\mathbb Z}^{d}\}=
   \{M\widetilde E' E'\beta+\widetilde E Es+c-\widetilde E Ec, \beta\in{\mathbb Z}^{d}\},$$
   where  $E'$,  $\widetilde E'$ are from~(\ref{fMEEM}). 
   On the other hand, 
   $\widetilde E( E \langle s\rangle)=\widetilde E  \langle q\rangle=
      \{\widetilde E M\beta+\widetilde E q+c-\widetilde E c, \beta\in{\mathbb Z}^{d}\}.$
      Since $q=Es+c-Ec-M\gamma,$ we see that 
         \begin{eqnarray}
          \widetilde E( E \langle s\rangle)= 
         \{\widetilde E M\beta+\widetilde E (Es+c-Ec-M\gamma)+c-\widetilde E c, \beta\in{\mathbb Z}^{d}\}=
         \nonumber
\\
         =\{ M (\widetilde E'\beta-\gamma)+\widetilde E Es+c-\widetilde E Ec, \beta\in{\mathbb Z}^{d}\}.     
          \nonumber         
         \end{eqnarray}
         Therefore, $\widetilde E( E \langle s\rangle)=\widetilde E E \langle s\rangle$
    $\Diamond.$
 	
Thus, Lemmas~\ref{lemCameron1} and~\ref{lemCameron2} are valid for the set ${\cal D}$. Hence, the set ${\cal D}$ can be split into the disjoint orbits, i.e. there exists a set $\Lambda\subset {\cal D}$ such that 
${\cal D}= \cup_{\langle s \rangle\in \Lambda}H\langle s \rangle.$ 
For convenience, redenote the elements of the set $\Lambda$ by $\langle s_{p,0}\rangle,$ where $p=0,\dots,\#\Lambda-1.$ 
As above, $H_{\langle s_{p,0}\rangle}$ is the stabilizer of $\langle s_{p,0}\rangle $. $\Gamma_{\langle s_{p,0}\rangle}$ is a complete set of representatives of the cosets of $H\backslash H_{\langle s_{p,0}\rangle}.$ 
The elements of the orbit $H\langle s_{p,0}\rangle$ we denote by $\langle s_{p,i}\rangle,$ where $i=0,\dots,\#\Gamma_{\langle s_{p,0}\rangle}-1.$ 
Let us fix index $p$. The matrices of the set 
$\Gamma_{\langle s_{p,0}\rangle}$ we denote by $E^{(i)}$ 
(or by $E_p^{(i)}$, $E_t^{(i)}$ when 
the matrices from the different sets $\Gamma_{\langle s_{p,0}\rangle}$ and $\Gamma_{\langle s_{t,0}\rangle}$ are used in a one equality), 
such that $E^{(i)} \langle s_{p,0}\rangle= \langle s_{p,i}\rangle, $ $E^{(0)}=I_d,$ $i=0,\dots,\#\Gamma_{\langle s_{p,0}\rangle}-1.$ 
The digit corresponding to the coset $\langle s_{p,i}\rangle$ we denote by $s_{p,i}.$

By Lemma~\ref{lemCameron2} for each 
$\langle s_{p,0}\rangle\in {\cal D} $ 
the symmetry group $H$ can be uniquely represented as
 follows $H=\Gamma_{\langle s_{p,0}\rangle} \times H_{\langle s_{p,0}\rangle},$ i.e. for each matrix $\widetilde E$ in $H$ there exist matrices $E\in \Gamma_{\langle s_{p,0}\rangle}$ and $F\in H_{\langle s_{p,0}\rangle}$ such that $\widetilde E=E F.$  
 The sets $\Gamma_{\langle s_{p,0}\rangle}, $ $H_{\langle s_{p,0}\rangle}$ can be considered as the "coordinate axes" of the symmetry group $H$. 
 (For each $p,$ 
 $p=0,\dots,\#\Lambda-1,$ these "coordinate axes" of $H$ can be different).

Now we 
point out some features of the choice of digits
and indicate how the  matrices from the  symmetry group $H$ act on the digits. 
  Let us fix a digit $s_{p,0}\in D(M)$ and a matrix $E^{(i)}\in\Gamma_{\langle s_{p,0}\rangle}.$ 
Since $E^{(i)} \langle s_{p,0}\rangle= \langle s_{p,i}\rangle,$ by~(\ref{fDigMain}) we get $E^{(i)} s_{p,0}+c-E^{(i)} c=M r^{E^{(i)}}_{p,0}+s_{p,i}$, where $r^{E^{(i)}}_{p,0}\in{\mathbb Z}^{d}.$ Let us rechoose the digits $s_{p,i}$ such that  $r^{E^{(i)}}_{p,0}=0,$ i.e.  
 \be
 E^{(i)} s_{p,0}+c-E^{(i)} c=:s_{p,i}, 
 \quad i=0,\dots,\#\Gamma_{\langle s_{p,0}\rangle}-1.
 \label{fDigE_i}
 \ee
 Throughout the paper we assume that the digits $s_{p,i}$ are chosen in that way.
  
 Suppose $F$ is a matrix from the stabilizer $H_{\langle s_{p,0}\rangle}$. Then $F\langle s_{p,0}\rangle=\langle s_{p,0}\rangle$. Hence by~(\ref{fDigMain}) 
   \be
  Fs_{p,0}=M r_{p,0}^F+ s_{p,0}+Fc-c,
  \label{fDigF}
  \ee
   where $r_{p,0}^F\in {\mathbb Z}^{d}$, $ r_{p,0}^F= M^{-1}(c-s_{p,0})-M^{-1}F(c-s_{p,0}).$
  
 For a matrix $K\in H$ there exist matrices 
 $E^{(i)}\in \Gamma_{\langle s_{p,0}\rangle}$ and 
 $F\in H_{\langle s_{p,0}\rangle}$ such that $K=E^{(i)} F$. Therefore, $K\langle s_{p,0}\rangle=
 E^{(i)} F\langle s_{p,0}\rangle=\langle s_{p,i}\rangle$. 
 Together with~(\ref{fDigE_i}) and~(\ref{fDigF}) we obtain that  
 	$$Ks_{p,0}=E^{(i)} Fs_{p,0}=E^{(i)} (M r_{p,0}^F+ s_{p,0}+Fc-c)=Mr^{K}_{p,0}+s_{p,i}+Kc-c,$$ 
 where $ r_{p,0}^K=M^{-1} E^{(i)} M r_{p,0}^F.$
  Analogously we can represent $Ks_{p,i}$. There exist matrices  
  $E^{(j(p,i,K))}\in \Gamma_{\langle s_{p,0}\rangle}$  and 
 $F\in H_{\langle s_{p,0}\rangle}$ such that $KE^{(i)}=E^{(j(p,i,K))} F$. Therefore, $K\langle s_{p,i}\rangle=
 E^{(j(p,i,K))} F \langle s_{p,0}\rangle=\langle s_{p,j(p,i,K)}\rangle$.
   Here the notation $j(p,\cdot,K)$ means the map from the set of indices $ \{0,\dots ,\#\Gamma_{\langle s_{p,0}\rangle}- 1\}$ to
itself and the index $j(p,i,K)$ is uniquely defined by the index $i$ and the matrix $K\in H$ for each $p.$ 
  Thus, together with~(\ref{fDigE_i}) and~(\ref{fDigF}) we obtain 
 \ba
  Ks_{p,i}=K E^{(i)} s_{p,0}+K c-K E^{(i)}  c=
  \hspace{5cm}
  \\  
   E^{(j)} M r_{p,0}^F + s_{p,j}+E^{(j)} F c-c +
   Kc -K E^{(i)} c=
   Mr^{K}_{p,i}+s_{p,j}+Kc-c,
     \nonumber
   \label{fDigK}
  \ea
   where $ r_{p,i}^K=M^{-1} E^{(j)} M r_{p,0}^F$ and $j=j(p,i,K). $

Let $H$ be a symmetry group on ${\mathbb Z}^{d}$ and 
    $t$ be an $H$-symmetric with respect to the center $c\in \frac 12 {\mathbb Z}^{d}$ trigonometric polynomial.
   Suppose that $H= H^{id}:=\{I_d, -I_d\}$, i.e. the simplest symmetry group. Then condition~(\ref{TrigSym}) is equivalent to

        \be
        t(\xi)=e^{2\pi i (2c,\xi)}t(-\xi) \quad\texttt{or}\quad h_k=h_{2c-k},\,\,\forall k\in{\mathbb Z}^{d}.
        \label{fSymMinus}
        \ee
    Such trigonometric polynomial $t(\xi)$ is called
    \emph{symmetric with respect to the point $c$} (or just point-symmetric). Easy to see that $H^{id}$ is a
    symmetry group with respect to any dilation matrix and Lemma~\ref{HanSym} for $H^{id}$ is always valid. 
    Note that if a symmetry group $H$ contains $H^{id}$ then  for any $H$-symmetric with respect to the point $c$ trigonometric polynomial $t$ we have

        $$D^{e_\ell} t({\bf0})=D^{e_\ell}\left(e^{2\pi i (2c,\xi)}t(-\xi)\right)\Big|_{\xi={\bf0}}=
        2\pi i 2(c)_j - D^{e_\ell}t({\bf0}),$$
    where for $\ell=1,\dots,d,$ the vectors $e_\ell\in{\mathbb Z}^{d}$ are the standard basis of ${\mathbb R}^d,$ $(e_\ell)_k=\delta_{k\ell},$ $k=1,\dots,d.$
    Therefore, $D^{e_\ell}t({\bf0})=2\pi i c^{e_\ell},\, \ell=1,\dots,d,$
    i.e. the phase for linear-phase moments of $t$
    must match with the point $c$. The maximal order of linear-phase moments for a point-symmetric
    trigonometric polynomial must be an even
    integer (see~\cite[Proposition 2]{Han11lpm}).

The following types of symmetry will be considered later.   The group
    $$H^{axis}:=\left\{ \texttt{diag}(u_1,\dots,u_d):
    u_j=\pm 1, j=1,\dots,d\right\}.$$
is called the axial symmetry group on ${\mathbb Z}^d$.
  The group
        $$H^{full}:=\left\{\pm I_2, \pm \left(\begin{matrix}
            -1 & 0\cr
            0 & 1\cr
        \end{matrix}\right), \pm \left(\begin{matrix}
            0 & 1\cr
            -1 & 0\cr
        \end{matrix}\right),
        \pm \left(\begin{matrix}
            0 & 1\cr
            1 & 0\cr
        \end{matrix}\right) \right\}$$
       is called the \emph{4-fold (or full) symmetry group on ${\mathbb Z}^2$.}

\begin{ex}  Let us demonstrate the above considerations on the choice of digits. Suppose $H=H^{full}$, $M=\left(\begin{matrix}
            2 & 0\cr
            0 & 2\cr
        \end{matrix}\right)$,
       $c=(0,0).$
       The set of digits is $D(M)=\{(0,0),(0,1),(1,0),(1,1)\}.$
       The set of cosets  is ${\cal D}:=\{\langle (0,0)\rangle,\langle (1,0)\rangle,\langle (0,1)\rangle,\langle (1,1)\rangle\}.$ The set of cosets is split into disjoint orbits. There are three orbits since 
       $$H\langle (0,0)\rangle=\{\langle (0,0)\rangle\},\,\, H\langle (1,1)\rangle=\{\langle (1,1)\rangle\},\,\, 
       H\langle (1,0)\rangle=\{\langle (1,0)\rangle,\langle (0,1)\rangle\}.$$
       Therefore, $\Lambda=\{\langle (0,0)\rangle,\langle (1,1)\rangle,\langle (1,0)\rangle\}.$
       Let us renumber the digits according with the  above considerations as follows
       $s_{0,0}=(0,0),$ $s_{1,0}=(1,1),$ $s_{2,0}=(1,0),$
       $s_{2,1}=(0,1).$
       For each digit $s_{p,0},$ $p=1,2,3,$ the sets $\Gamma_{\langle s_{p,0}\rangle}, $ $H_{\langle s_{p,0}\rangle}$ can be defined.
       Namely, 
       $$H_{\langle s_{0,0}\rangle}=H^{full}, \,\,
       \Gamma_{\langle s_{0,0}\rangle}=\{I_2\},\,\,
       H_{\langle s_{1,0}\rangle}=H^{full}, \,\,
       \Gamma_{\langle s_{1,0}\rangle}=\{I_2\},\,\,$$      
       $$H_{\langle s_{2,0}\rangle}=
       \left\{\pm I_2, \pm \left(\begin{matrix}
            1 & 0\cr
            0 & -1\cr
        \end{matrix}\right)\right\}, \,\,
       \Gamma_{\langle s_{2,0}\rangle}= 
       \left\{ I_2, \left(\begin{matrix}
            0 & 1\cr
            1 & 0\cr
        \end{matrix}\right)\right\}.$$
        
\end{ex}

\section{Preliminary results}

The general scheme for the construction of compactly supported MRA-based dual wavelet systems 
    was developed in~\cite{RS},
    \cite{RS2} (the Mixed Extension Principle).
    To construct the pair of such wavelet systems 
    one
    starts with two compactly supported refinable functions
    $\phi, \widetilde\phi$, $\widehat\phi({\bf0})=1$, $\widehat{\widetilde\phi}({\bf0})=1,$
    (or its masks $m_0, \widetilde m_0$, respectively, which are trigonometric
    polynomials).
    Then one finds trigonometric
    polynomials $m_{\nu}, \widetilde m_{\nu}$, $\nu = 1,\dots, r$,
    $r\ge m-1$, called the wavelet masks, such that the following polyphase matrices 
        \ban
        \cal M:=\left(\begin{array}{ccc}
        \mu_{00}&\hdots & \mu_{0, m-1}\cr
        \vdots & \ddots & \vdots \cr
        \mu_{r,0}& \hdots & \mu_{r,m-1}
        \end{array}\right),\,
        \cal \widetilde M:=\left(\begin{array}{ccc}
        \widetilde{\mu}_{00}&\hdots & \widetilde{\mu}_{0, m-1}\cr
        \vdots & \ddots & \vdots \cr
        \widetilde{\mu}_{r,0}& \hdots & \widetilde{\mu}_{r,m-1}
        \end{array}\right)
        \ean
    satisfy
        \be
        {\cal M}^*{\cal\widetilde M}= I_m
        \label{calM}.
        \ee
    Here $\mu_{\nu k},$ $\widetilde\mu_{\nu k}$ $k=0,\dots,m-1,$ are the polyphase components of the
    wavelet masks $m_{\nu},$ $\widetilde m_{\nu}$ for all $\nu=0,\dots,r,$ $r\ge m-1.$

The wavelet functions $\psi^{(\nu)}$, $\widetilde\psi^{(\nu)},$
    $\nu = 1,\dots, r$, $r\ge m-1,$ are defined via its Fourier transform
        \be \widehat{\psi^{(\nu)}}(\xi)=
        m_\nu(M^{*-1}\xi)\widehat\varphi(M^{*-1}\xi), \quad
        \widehat{\widetilde\psi^{(\nu)}}(\xi)=
        \widetilde m_\nu(M^{*-1}\xi)\widehat{\widetilde\varphi}(M^{*-1}\xi).
        \label{PsiDef}
        \ee
If the wavelet functions  $\psi^{(\nu)}, \widetilde\psi^{(\nu)}$,
    $\nu=1,\dots, r$, $r\ge m-1,$
    are constructed as above, then
    the set of the functions $\{\psi_{jk}^{(\nu)}\}$, $\{\widetilde\psi_{jk}^{(\nu)}\}$
    is said to be a
    \emph{compactly supported MRA-based dual wavelet system} generated by the
    refinable functions $\varphi, \widetilde\varphi$
    (or their masks $m_0, \widetilde m_0$).

For the  wavelet functions $\widetilde\psi^{(\nu)}$, $\nu=1,\dots, r,$ $r\ge m-1,$ defined by~(\ref{PsiDef}),
   the $VM^{n}$ property for the corresponding
    wavelet system is equivalent to the fact that the wavelet masks $\widetilde m_\nu$, $\nu=1,\dots, r$
    have vanishing moments up to order $n$.

Suppose that some refinable functions $\varphi, \widetilde\varphi\in L_2({\mathbb R}^d)$ generate a compactly supported MRA-based dual wavelet system $\{\psi_{jk}^{(\nu)}\}$, $\{\widetilde\psi_{jk}^{(\nu)}\},$
    $\nu=1,\dots, r$, $r\ge m-1.$
    A necessary (see~\cite[Theorem 1]{Sk1}) and
    sufficient (see~\cite[Theorems 2.2, 2.3]{Han1}) condition
    for the wavelet system $\{\psi_{jk}^{(\nu)}\}$, $\{\widetilde\psi_{jk}^{(\nu)}\}$
    to be a pair of dual wavelet frames in $L_2({\mathbb R}^d)$ is that
    each wavelet system $\{\psi_{jk}^{(\nu)}\}$ and $\{\widetilde\psi_{jk}^{(\nu)}\}$
    has vanishing moments at least of order $1$.
     For instance, in the multidimensional case the explicit method for the construction
    of compactly supported dual wavelet frames with vanishing moments was suggested in~\cite{Sk1}.
    But it is a difficult task to provide various types of symmetry for the wavelet frames.
    
    For the special case of interpolatory refinable masks the scheme can be simplified. 
    Recall that if $m_0$ is interpolatory then $\mu_{00}\equiv \frac 1{\sqrt m}.$ Moreover, if $m_0$ has sum rule of order $n\in{\mathbb N}$ then by Lemma~\ref{lemDynSk} we have that $\lambda_\beta=\delta_{{\bf0}\beta}$
and conditions~(\ref{fPolySRInterp}) are valid.
The dual mask can be defined in the following form via its polyphase components
\be(\widetilde\mu_{00}, \widetilde\mu_{01},\dots,\widetilde\mu_{0,m-1})=\left(\sqrt m \left(1-\sum\limits_{k=1}^{m-1}|\mu_{0k}|^2 \right), \mu_{01},\dots,\mu_{0,m-1}\right).
\label{fWM0Def}
\ee
Not hard to check that with this definition we provide
$
\sum_{k=0}^{m-1}\mu_{0k}\overline{\widetilde\mu_{0k}}\equiv 1.
$
Also, by~(\ref{fPolySRInterp}) $\mu_{0k}({\bf0})=\frac 1{\sqrt{m}},$ $k=0,\dots,m-1.$ Then  $\widetilde\mu_{00}({\bf0})=\frac 1{\sqrt{m}}$ and by Lemma~\ref{lemDynSk} we can conclude that $\widetilde m_0$  has sum rule at least of
 order $1$ and $\widetilde \lambda_{\beta}=\delta_{{\bf0}\beta}.$
 In fact,  the mask $\widetilde m_0$ has sum rule of order $n$ (see~\cite[Proposition 6]{Sk6}).


The order of sum rule of refinable masks is connected with the order of vanishing moments of wavelet masks. The order of vanishing moments is connected with the approximation order of the corresponding wavelet system (see, e.g.~\cite{Sk1}).
The method providing the vanishing moments for the wavelet masks is given by

\begin{lem} (\cite[Theorem 10]{KST})
Let  $\phi, \widetilde\phi\in S'$ be refinable functions with compact support. 
Let $m_0, \widetilde m_0$ be their refinable masks. Suppose $m_0$ is interpolatory and has sum rule of order $n$. Mask $\widetilde m_0$ is defined by~(\ref{fWM0Def}).
    Then there exist a MRA-based dual wavelet system
    $\{\psi_{jk}^{(\nu)}\}$, $\{\widetilde\psi_{jk}^{(\nu)}\}$, $\nu=1,\dots, m-1,$
    such that the wavelet system $\{\widetilde\psi^{(\nu)}_{jk}\}$
    has the $VM^{n}$ property.
   \label{lemSkVM}
\end{lem}

 Note that if $\phi, \widetilde\phi$ are in $L_2({\mathbb R}^d),$
 then the  MRA-based dual wavelet system in Lemma~\ref{lemSkVM} is a dual wavelet frame.
The technique for the extension of the polyphase matrices
   can be  realised as follows. Let us define the part of the row of the polyphase components of $m_0$ by $P=(\mu_{01},\dots,\mu_{0,m-1}).$ 
 Then  the explicit formulas  for the  matrix extension are given by 
\be
   {\cal N}=\left( \begin{array}{c|c}
      \frac 1 {\sqrt m} & P \\
      \hline
      &  \\
     -U P^{*} & \sqrt m(U-U P^{*}  P)
    \end{array}\right),
\quad 
{\cal \widetilde N}=\left(
\begin{array}{c|c}
     \sqrt m (1-PP^*)& P \\
      \hline
      &   \\
      -\widetilde UP^{*} &  \frac 1 {\sqrt m} \widetilde U
    \end{array} \right)
\label{fNNblockDual} 
 \ee
    Here $U,$ $\widetilde U$ are $(m-1)\times (m-1)$ matrices consisting of trigonometric polynomials such that
 $U\widetilde U^*\equiv I_{m-1}.$ It follows easily that  ${\cal N}\widetilde {\cal N}^*\equiv I_{m}.$ The matrices $U,$ $\widetilde U$, for example, can be taken  as follows $U=\widetilde U=I_{m-1}.$
 The next section is devoted to  the construction of wavelets using the above technique such that the resulting wavelet systems have some symmetry properties.

\section{Multivariate $H$-symmetric dual wavelet frames}

Firstly, we state two simple lemmas.

\begin{lem}
        Let $t'(\xi)$ be a trigonometric polynomial, $n\in{\mathbb N},$
            $t(\xi):=e^{2\pi i (a,\xi)}t'(\xi),$ where $a\in{\mathbb R}^d.$
        Then
            $D^{\beta}t'(\xi)\Big|_{\xi={\bf0}}=(2\pi i)^{[\beta]}\kappa'_{\beta}, $
             $\forall \beta\in\Delta_n,$
        if and only if
            $$D^{\beta}t(\xi)\Big|_{\xi={\bf0}}=(2\pi i)^{[\beta]}\sum\limits_{\alpha\in\square_{\beta}}
            \kappa'_{\alpha} \binom{\beta}{\alpha} a^{\beta-\alpha}, \quad  \forall \beta\in\Delta_n.$$

    \label{LambdaLem}
\end{lem}
      Let us fix a symmetry group $H$ on ${\mathbb Z}^{d}$, an appropriate symmetry center $c\in{\mathbb R}^d$ and an appropriate matrix dilation $M$ with the set of digits $D(M).$ Now we reformulate the 
    $H$-symmetry condition for a trigonometric polynomial in terms of its polyphase components.
    The polyphase components of trigonometric polynomial $t$ we enumerate as the corresponding digits: $\tau_{p,i}(\xi),$
   $i=0,\dots,\#\Gamma_{\langle s_{p,0}\rangle}-1,$
     $p=0,\dots,\#\Lambda-1$.

\begin{lem}
    A trigonometric polynomial $t$ is $H$-symmetric with respect to
    the appropriate center $c\in{\mathbb R}^d$
    if and only if
    for each $p\in\{0,\dots,\#\Lambda-1\}$ its polyphase components $\tau_{p,i},$  
    $i=0,\dots,\#\Gamma_{\langle s_{p,0}\rangle}-1,$
    satisfy 
        \be
        \tau_{p,i}((M^{-1} K M)^*\xi)=e^{2\pi i (r^K_{p,i},\xi)} \tau_{p,j}(\xi), \quad \forall K\in H
        \label{fLemSymPoly}
        \ee
        where index $j=j(p,i,K)$ is defined as above in~(\ref{fDigK}), 
         vector $ r_{p,i}^K=M^{-1} E^{(j)} M r_{p,0}^{F},$ $r_{p,0}^{F}$ is from~(\ref{fDigF}), the matrices  
  $E^{(j)}\in \Gamma_{\langle s_{p,0}\rangle}$  and 
 $F\in H_{\langle s_{p,0}\rangle}$ are such that $KE^{(i)}=E^{(j)} F$,   $E^{(i)}\in \Gamma_{\langle s_{p,0}\rangle}.$ 

\label{lemSymPoly}
\end{lem}

{\bf Proof.} Suppose $t$ is an $H$-symmetric with respect to the center $c$ trigonometric polynomial.          
Then by~(\ref{PR}) and~(\ref{fDigK}) we obtain
                $$ t(K^*\xi)= \frac1{\sqrt m}\sum\limits_{p=0}^{\#\Lambda-1} \sum\limits_{i=0}^{\#\Gamma_{\langle s_{p,0}\rangle}-1} 
                e^{2\pi i(K s_{p,i},\xi)}
                \tau_{p,i}(M^*K^*\xi)=$$
                $$\frac1{\sqrt m}
                e^{2\pi i (Kc-c,\xi)}
                \sum\limits_{p=0}^{\#\Lambda-1}
                 \sum\limits_{i=0}^{\#\Gamma_
                 {\langle s_{p,0}\rangle}-1} 
                 e^{2\pi i(s_{p,j},\xi)}
                 \tau_{p,i}(M^*K^*\xi) 
                 e^{2\pi i(M r^K_{p,i},\xi)},$$
                $$e^{2\pi i (Kc-c,\xi)} t(\xi)=\frac1{\sqrt m}e^{2\pi i (Kc-c,\xi)}\sum\limits_{p=0}^{\#\Lambda-1} \sum\limits_{i=0}^{\#\Gamma_{\langle s_{p,0}\rangle}-1} 
                e^{2\pi i(s_{p,j},\xi)}
                \tau_{p,j}(M^*\xi),              $$
                where $K\in H$, the index $j=j(p,i,K)$ is defined  as above, i.e
               for each $p\in\{0,\dots,\#\Lambda-1\}$
                there exist matrices  
  $E^{(j(p,i,K))}\in \Gamma_{\langle s_{p,0}\rangle}$  and 
 $F\in H_{\langle s_{p,0}\rangle}$ such that $KE^{(i)}=E^{(j(p,i,K))} F$. 
  Since the polyphase representation is unique with respect to the chosen digits, we get~(\ref{fLemSymPoly}).
The converse statement is obvious.
$\Diamond$

According with the above considerations we can rewrite condition~(\ref{fLemSymPoly}) in the following form.
If $F\in H_{\langle s_{p,0}\rangle}$ then
        \be
        \tau_{p,0}(\xi)=
        e^{2\pi i (r^F_{p,0},\xi)}\tau_{p,0}((M^{-1} F M)^*\xi),
        \label{fPolyF}
        \ee
where $r^F_{p,0}$ is from~(\ref{fDigF}).  Thus, $\tau_{p,0}(\xi)$ is $M^{-1}H_{\langle s_{p,0}\rangle}M-$symmetric with respect to the center $M^{-1}(c-s_{p,0}).$
If $E^{(i)}\in\Gamma_{\langle s_{p,0}\rangle}$, 
 $i=0,\dots,\#\Gamma_{\langle s_{p,0}\rangle}-1,$ then 
        \be
        \tau_{p,i}(\xi)=\tau_{p,0}((M^{-1} E^{(i)} M)^*\xi).
        \label{fPolyE_i}
        \ee
        Note that $\tau_{p,i}$ is  
        $M^{-1}H_{\langle s_{p,i}\rangle}M-$symmetric with respect to the center $M^{-1}(c-s_{p,i}).$
Indeed, since the group $H_{\langle s_{p,i}\rangle}$ is conjugate to $
H_{\langle s_{p,0}\rangle}$, i.e. $H_{\langle s_{p,i}\rangle}=E^{(i)} H_{\langle s_{p,0}\rangle} (E^{(i)})^{-1}$,
we see that for any $K\in H_{\langle s_{p,i}\rangle}$ there exists $F\in H_{\langle s_{p,0}\rangle}$
such that $K=E^{(i)} F (E^{(i)})^{-1}.$
Thus, we need to show that $\tau_{p,i}((M^{-1} K M)^*\xi)e^{2\pi i (r^K_{p,i},\xi)}=\tau_{p,i}(\xi),$
where
$$r^K_{p,i}=M^{-1}(c-s_{p,i})-(M^{-1}K M) M^{-1}(c-s_{p,i}).$$

 Indeed, 
by~(\ref{fPolyF}) and~(\ref{fPolyE_i}) we obtain
	$$\tau_{p,i}((M^{-1} K M)^*\xi)e^{2\pi i (r^K_{p,i},\xi)}
	=\tau_{p,0}(M^* F^* E^{(i)*} M^{*-1}\xi)e^{2\pi i (r^K_{p,i},\xi)}=$$
	$$\tau_{p,0}(M^* E^{(i)*} M^{*-1}\xi)
	e^{-2\pi i (r_{p,0}^F,M^*E^{(i)*} M^{*-1}\xi)}e^{2\pi i (r^K_{p,i},\xi)}
	=\tau_{p,i}(\xi).$$
The next to last equality is valid due to 
$r^K_{p,i}=M^{-1} E^{(i)} Mr_{p,0}^F.$ 


\subsection{Construction of symmetric refinable masks}

   Now we suggest a simple
    algorithm for the construction of interpolatory refinable masks that are
    $H$-symmetric with respect to
    the origin ($c={\bf0}$) and
    have an arbitrary order of sum rule.
    Assume that $s_{0,0}={\bf0},$ then $H_{\langle s_{0,0}\rangle}=H$, 
    $\Gamma_{\langle s_{0,0}\rangle}=I_d.$
    Conditions~(\ref{fDigE_i}),~(\ref{fDigF}),~(\ref{fDigK}) 
    in this case become
    $ E^{(i)} s_{p,0}=s_{p,i}$,
    $ Fs_{p,0}=M r_{p,0}^F+ s_{p,0},$
   $Ks_{p,i}=Mr^{K}_{p,i}+s_{p,j},$ where $j=j(p,i,K).$
Also, since in the interpolatory case $\lambda_{\gamma}=\delta_{{\bf0}\gamma}$ then condition~(\ref{fPolySumRule}) is equivalent to condition~(\ref{fPolySRInterp}).

\begin{theo}
    Let $H$ be a symmetry group on ${\mathbb Z}^{d}$, $M$ be an appropriate matrix dilation, $n\in{\mathbb N}$.
    Then there exists an interpolatory mask $m_0$
     that is $H$-symmetric with
    respect to the origin and
    has the order of sum rule equal to $n.$
    This mask $m_0$ can be represented by 
        \ba
        m_0(\xi)= \frac 1m + \frac1{\sqrt m}
        \sum\limits_{p=1}^{\#\Lambda-1} \nu_p(\xi).    
        \label{fm_0HSym}
        \ea
    Here for each $p=1,\dots,\#\Lambda-1$     
    	$$
    	\nu_p(\xi)     =
    	\sum\limits_{i=0}^{\#\Gamma_{\langle s_{p,0}\rangle}-1}      
        \frac { e^{2\pi i(s_{p,i},\xi)}} 
        {\# H_{\langle s_{p,0}\rangle}}
        \sum\limits_{F\in H_{\langle s_{p,0}\rangle}}
        G_{0,p,0}(M^* F^* E^{(i)*}\xi)
        e^{2\pi i (r_{p,0}^F,M^* E^{(i)*}\xi)},
    	$$
   where $G_{0,p,0}$ are trigonometric polynomials such that~(\ref{fG0p0}) is valid, 
    $E^{(i)}\in \Gamma_{\langle s_{p,0}\rangle}.$

    \label{theoMaskH}
\end{theo}

{\bf Proof.}
Let us construct the polyphase components  $\mu_{0,p,i}(\xi),$ 
 $i=0,\dots,\#\Gamma_{\langle s_{p,0}\rangle}-1,$
     $p=0,\dots,\#\Lambda-1$ such that
they satisfy conditions~(\ref{fPolySRInterp}) and~(\ref{fLemSymPoly}). Therefore,
we obtain the required mask $m_0$ by~(\ref{PR}).

For $p=0$ we can set $\mu_{0,0,0}$ equal to $\frac1{\sqrt m}.$

Let us fix $p\neq0$ and
construct the polyphase component $\mu_{0,p,0}(\xi)$. It should satisfy conditions~(\ref{fPolySRInterp}) and~(\ref{fPolyF}).
First of all, define the trigonometric polynomial $G_{0,p,0}(\xi)$ such that 
	\be
	D^{\beta}G_{0,p,0}({\bf0})=
	\frac {(2\pi i)^{[\beta]}} {\sqrt m}
        	(-M^{-1}s_{p,0})^{\beta}
        	\quad \forall \beta\in\Delta_n.
	\label{fG0p0}        
        \ee
Thus, $\mu_{0,p,0}$ can be obtained as follows
        \be
        \mu_{0,p,0}(\xi)=\frac 1 {\# H_{\langle s_{p,0}\rangle}}
        \sum\limits_{F\in H_{\langle s_{p,0}\rangle}}
        G_{0,p,0}((M^{-1} F M)^*\xi)e^{2\pi i (r_{p,0}^F,\xi)}.
        \label{fMu0p0}
        \ee
Let us check that $\mu_{0,p,0}(\xi)$ is $M^{-1}H_{\langle s_{p,0}\rangle}M-$symmetric with respect to the center $-M^{-1}s_{p,0}$ and conditions~(\ref{fPolySRInterp}) are valid.
Suppose $\widetilde F\in H_{\langle s_{p,0}\rangle}$. Then the  $M^{-1}H_{\langle s_{p,0}\rangle}M$-symmetry of $\mu_{0,p,0}$
 follows from
$$\mu_{0,p,0}((M^{-1} \widetilde F M)^*\xi)=
	\frac 1 {\# H_{\langle s_{p,0}\rangle}}
        \sum\limits_{F\in H_{\langle s_{p,0}\rangle}}
        G_{0,p,0}((M^{-1} \widetilde F F M)^*\xi)e^{2\pi i (M^{-1} \widetilde F Mr_{p,0}^F),\xi)}=$$
        $$\frac 1 {\# H_{\langle s_{p,0}\rangle}}
        \sum\limits_{F\in H_{\langle s_{p,0}\rangle}}
        G_{0,p,0}((M^{-1} \widetilde F F M)^*\xi)
        e^{2\pi i (r_{p,0}^{\widetilde FF}-r_{p,0}^{\widetilde F},\xi)}=
        \mu_{0,p,0}(\xi)e^{-2\pi i (r_{p,0}^{\widetilde F},\xi)},$$
since $M^{-1} \widetilde F Mr_{p,0}^F=r_{p,0}^{\widetilde FF}-r_{p,0}^{\widetilde F}.$       
By Lemma~\ref{LambdaLem}, condition~(\ref{fG0p0}) is equivalent to
$$D^{\beta}\left(G_{0,p,0}(\xi) 
e^{2\pi i (M^{-1}s_{p,0},\xi)}\right)
\Big|_{\xi={\bf0}}=\frac {(2\pi i)^{[\beta]}} {\sqrt m} \delta_{{\bf0}\beta}
\quad \forall \beta\in\Delta_n.$$
Due to the higher chain rule for the linear change of variables
 $\xi:=(M^{-1} F M)^*\xi,$ $F\in H_{\langle s_{p,0}\rangle}$ we have 
$$D^{\beta}\left(G_{0,p,0}((M^{-1} F M)^*\xi) 
e^{2\pi i (M^{-1}Fs_{p,0},\xi)}\right)
\Big|_{\xi={\bf0}}=\frac {(2\pi i)^{[\beta]}} {\sqrt m}  \delta_{{\bf0}\beta}
\quad \forall \beta\in\Delta_n.$$
Then by Lemma~\ref{LambdaLem} and~(\ref{fDigF}) we conclude that $\forall F\in H_{\langle s_{p,0}\rangle}$ 
$$D^{\beta}\left(G_{0,p,0}((M^{-1} F M)^*\xi) 
e^{2\pi i (r_{p,0}^F,\xi)}\right)
\Big|_{\xi={\bf0}}=\frac {(2\pi i)^{[\beta]}} {\sqrt m}
        (-M^{-1}s_{p,0})^{\beta} \quad \forall \beta\in\Delta_n.$$
        Hence, condition~(\ref{fPolySRInterp}) for $\mu_{0,p,0}$ is valid.

Next, define $\mu_{0,p,i}(\xi)$ as follows
$\mu_{0,p,i}(\xi)=\mu_{0,p,0}((M^{-1} E^{(i)} M)^*\xi),$ where $E^{(i)}\in \Gamma_{\langle s_{p,0}\rangle}.$ 
Then the polyphase component  $\mu_{0,p,i}(\xi)$ is $M^{-1}H_{\langle s_{p,i}\rangle}M-$symmetric with respect to the center $-M^{-1}s_{p,i}.$

Show that condition~(\ref{fPolySRInterp}) is valid for $\mu_{0,p,i}(\xi).$
We have already proved that 
$$D^{\beta}\left(\mu_{0,p,0}(\xi) 
e^{2\pi i (M^{-1}s_{p,0},\xi)}\right)
\Big|_{\xi={\bf0}}=\frac {(2\pi i)^{[\beta]}} {\sqrt m} \delta_{{\bf0}\beta}
 \quad \forall \beta\in\Delta_n.$$
Hence, by the higher chain rule for the linear change of variables
 $\xi:=(M^{-1} E_p^{(i)} M)^*\xi$, $E_p^{(i)}\in\Gamma_{\langle s_{p,0}\rangle}$, we get 
$$D^{\beta}\left(\mu_{0,p,0}((M^{-1} E^{(i)} M)^*\xi) 
e^{2\pi i (M^{-1} E^{(i)}s_{p,0},\xi)}\right)
\Big|_{\xi={\bf0}}=\frac {(2\pi i)^{[\beta]}} {\sqrt m} \delta_{{\bf0}\beta}
 \quad \forall \beta\in\Delta_n.$$
It remains to note that $E^{(i)}s_{p,0}=s_{p,i}$ and then by Lemma~\ref{LambdaLem}
$$D^{\beta}\mu_{0,p,i}({\bf0})=\frac {(2\pi i)^{[\beta]}} {\sqrt m}
       (-M^{-1}s_{p,i})^{\beta}  \quad \forall \beta\in\Delta_n.$$
        
Hence, we construct the polyphase components $\mu_{0,p,i}(\xi),$ 
 $i=0,\dots,\#\Gamma_{\langle s_{p,0}\rangle}-1,$
     $p=0,\dots,\#\Lambda-1$  such that they
    satisfy symmetry condition~(\ref{fLemSymPoly}) in Lemma~\ref{lemSymPoly}
    and condition~(\ref{fPolySRInterp}) with $\lambda_{\gamma}=\delta_{{\bf0}\gamma}$.
    To get the required mask $m_0$ it remains to use formula~(\ref{PR}) together with~(\ref{fMu0p0}).
    $\Diamond$

\begin{rem}
    The construction of $H$-symmetric with
    respect to the point $c$ masks with real coefficients and with  sum rule of order $n$ 
    can be done as in
    Theorem~\ref{theoMaskH} using~(\ref{fm_0HSym}), where $G_{0,p,0}$ 
    are the same trigonometric polynomial  but with real coefficients.
    Also, the minimal number of non-zero coefficients for
    masks from Theorem~\ref{theoMaskH} can be achieved if  
    we provide the minimality of non-zero coefficients  for 
    $\mu_{0,p,0},$ $p=0,\dots,\#\Lambda-1$.
    \label{rRealCoef}
\end{rem}

\subsection{Construction of symmetric wavelet frames}

For the construction of wavelet frames with the symmetry properties we need an $H$-symmetric dual refinable mask $\widetilde m_0.$ It can be  defined via its polyphase components by~(\ref{fWM0Def}). As it was mentioned  above $\widetilde m_0$ has the same order of sum rule as $m_0.$ If $m_0$ is $H$-symmetric then $\widetilde m_0$ is also $H$-symmetric. Indeed, by~(\ref{fLemSymPoly}) it is enough to show that $\widetilde\mu_{00}$ is $H$-symmetric. But this is true since $m_0$ is $H$-symmetric.

Next, we consider the basic matrix extension as in~(\ref{fNNblockDual}) with $U=\widetilde U=I_{m-1}.$ Namely,
\be
   {\cal N}=\left( \begin{array}{c|c}
      \frac 1 {\sqrt m} & P \\
      \hline
      &  \\
     - P^{*} & \sqrt m(I_{m-1}- P^{*}  P)
    \end{array}\right),
\quad 
{\cal \widetilde N}=\left(
\begin{array}{c|c}
     \sqrt m (1-PP^*)& P \\
      \hline
      &   \\
      -P^{*} &  \frac 1 {\sqrt m} I_{m-1}
    \end{array} \right)
\label{fNNblockDualIdNoS} 
 \ee
 
 \begin{theo} 
  Let $H$ be a symmetry group on ${\mathbb Z}^{d}$, 
  $M$ be an appropriate matrix dilation, $n\in{\mathbb N}$
and let $m_0$ and $\widetilde m_0$ be $H$-symmetric
    with respect to the origin masks such that
     $m_0$
    has sum rule of order $n$,
      $\widetilde m_0$ is constructed by~(\ref{fWM0Def}).
    Then wavelet masks  $m_{(p,i)}$ and $\widetilde m_{(p,i)},$
    $i=0,\dots,\#\Gamma_{\langle s_{p,0}\rangle}-1,$
     $p=1,\dots,\#\Lambda-1$,
    constructed using the matrix extension~(\ref{fNNblockDualIdNoS})
  have the $VM^n$ property and the symmetry properties. Namely, for a fixed number $p=1,\dots,\#\Lambda-1$,
   $m_{(p,i)}$ and $\widetilde m_{(p,i)}$ are
    $H_{\langle s_{p,i}\rangle}$-symmetric with respect to
    the center $s_{p,i}$; 
      $m_{(p,i)}=m_{(p,0)}(E^{(i)*}\xi),$
        $\widetilde m_{(p,i)}= 
      \widetilde m_{(p,0)}(E^{(i)*}\xi),$ where 
    $E^{(i)}\in\Gamma_{\langle s_{p,0}\rangle},$
    $i=0,\dots,\#\Gamma_{\langle s_{p,0}\rangle}-1$. 
     If $\phi,$ $\widetilde\phi \in L_2({\mathbb R}^d)$, then the corresponding 
     wavelet system is a dual wavelet frame.
\label{theoWaveNoSym}
\end{theo}

{\bf Proof.} 
Let $T_p$ denote the row of the polyphase components defined by 
$T_p=(\mu_{0,p,0},...,
\mu_{0,p,N_p-1}),$
 where $N_p=\#\Gamma_{\langle s_{p,0}\rangle},$ $p=1,...,\#\Lambda-1.$
  Let $P=(T_1,\dots,T_{\#\Lambda-1})$.
Consider the submatrix $I_m-P^*  P$ in~(\ref{fNNblockDualIdNoS})

    $$I_m-P^{*}  P=
    \left( \begin{array}{cccc}
      I_{N_1-1}-
      T_1^*  T_1 & 
      	-T_1^* T_2 & \dots & 
      		-T_1^* T_{\#\Lambda-1}\\
      -T_2^* T_1 &
      	I_{N_2-1}-T_2^* T_2 & \dots & 
      		-T_2^* T_{\#\Lambda-1} \\
      		\vdots & & \ddots & \vdots \\
      	-T_{\#\Lambda-1}^*  T_1 &	\dots &\dots &
      	I_{N_{(\#\Lambda-1)}-1}-T_{\#\Lambda-1}^*
      	 T_{\#\Lambda-1}
    \end{array}\right).$$
    The rows of the submatrix are numbered by the double indices $(p,i)$ corresponding to the indices of the polyphase components $\overline{\mu_{0,p,i}}$ in the columns $T_p^*,$ $p=1,...,\#\Lambda-1.$
    Consider the row $(p,i).$
    The corresponding wavelet mask for this row is defined by
    $m_{(p,i)}=e^{2\pi i (s_{p,i},\xi)}-\overline{\mu_{0,p,i}(M^*\xi)}m_0(\xi).$
    The mask $m_{(p,0)}$ is $H_{\langle s_{p,0}\rangle}$-symmetric with respect to the center $s_{p,0}.$
    Indeed, due to~(\ref{fDigF}) and~(\ref{fPolyF}) for all 
      $F\in H_{\langle s_{p,0}\rangle}$ we have
  
    $$ m_{(p,0)}(F^*\xi)= 
    e^{2\pi i (F s_{p,0},\xi)} - \overline{\mu_{0,p,0}(M^*F^*\xi)} m_0(F^*\xi)=
    m_{(p,0)}(\xi)e^{2\pi i (Fs_{p,0}-s_{p,0},\xi)}.
    $$
    Moreover, by~(\ref{fPolyE_i}) we have
   $m_{(p,i)}(\xi)=m_{(p,0)}(E^{(i)*}\xi).$
    Similarly, it can be proved that 
    $m_{(p,i)}$ is $H_{\langle s_{p,i}\rangle}$-symmetric with respect to the center $s_{p,i}.$ Namely, for all
    $\widetilde F \in H_{\langle s_{p,i}\rangle}$
     $$ m_{(p,i)}(\widetilde F^*\xi)=m_{(p,i)}(\xi)e^{2\pi i (\widetilde Fs_{p,i}-s_{p,i},\xi)}. $$
     Easy to check that the same symmetric properties are valid for the dual masks, since 
     $$\widetilde m_{(p,i)}=\frac {e^{2\pi i (s_{p,i},\xi)}}{m}-\frac {\overline{\mu_{0,p,i}(M^*\xi)}}{\sqrt m}.$$
     The $VM^n$ property for the  wavelet masks is provided by Lemma~\ref{lemSkVM}.
     $\Diamond$
     
    Thus, the matrix extension as in~(\ref{fNNblockDualIdNoS}) leads to the wavelet masks the are mutually symmetric, i.e. some wavelet masks are reflected or rotated copies of the others.  Proposition 2.1 in~\cite{Han3} can be used in order to compute the symmetry centers of wavelet functions. Notice that if $\Gamma_{\langle s_{p,0}\rangle}=\{I_d\}$ for all 
    $p=0,\dots,\#\Lambda-1$, then all wavelet masks constructed by
    Theorem~\ref{theoWaveNoSym}  are $H$-symmetric.
    
    \begin{rem}  
    All wavelet functions constructed in Theorem~\ref{theoWaveNoSym} have the $VM^n$ property.
    It is known that the dual wavelet frame provides approximation order $n$ if the dual wavelet functions $\{\widetilde \psi^{(\nu)}\}$ have the $VM^n$ property (see Theorem 4 in~\cite{Sk1}). So, to obtain dual wavelet frames we just  need to require the vanishing moments at least of order 1 for the primal wavelet functions $\{ \psi^{(\nu)}\}$ (as noted above, this is a necessary and sufficient condition for a wavelet system to be a frame).
   Thus, it is enough to provide the 
   sum rule at least of order 1 for the dual refinable mask $\widetilde m_0$.  The general method for the construction of dual refinable masks with desirable order of sum rule was suggested by B. Han in~\cite{HanAnalSmooth2}.
  Here we slightly modify this method by adding the $H$-symmetry conditions. 
  Let us construct $\widetilde m_0$ using Theorem~\ref{theoMaskH} with $n=1,$ $\widetilde{\mu}_{0k}$   be its polyphase components, $k=0,\dots,m-1$. Now we modify the mask $\widetilde m_0$ by replacing its first polyphase component $\widetilde{\mu}_{00}=\frac 1{\sqrt m}$ with
  $\widetilde{\mu}_{00}= \sqrt m ( 1-\sum_{k=1}^{m-1}
\overline{\mu_{0k}} \widetilde{\mu}_{0k})$,
where $\mu_{0k}$ are the polyphase components
of the mask $m_0$. 
Easy to see that the modified mask $\widetilde m_0$ satisfies the following 
    conditions:
    $\widetilde m_0$ is $H$-symmetric with respect to the origin,
    $\widetilde m_0$ has sum rule at least of order 1 and
    $
\sum_{k=0}^{m-1}\mu_{0k}\overline{\widetilde\mu_{0k}}\equiv 1.
$
Analyzing the proof of Theorem~\ref{theoWaveNoSym} we can restate it for the modified dual mask $\widetilde m_0$. The matrix extension in this case is as follows
\be
   {\cal N}=\left( \begin{array}{c|c}
      \frac 1 {\sqrt m} & P \\
      \hline
      &  \\
     - \widetilde P^{*} & \sqrt m(I_{m-1}- \widetilde P^{*}  P)
    \end{array}\right),
\quad 
{\cal \widetilde N}=\left(
\begin{array}{c|c}
     \sqrt m (1- \widetilde P  P^*)& \widetilde P \\
      \hline
      &   \\
      -P^{*} &  \frac 1 {\sqrt m} I_{m-1}
    \end{array} \right),
\label{fNNblockDualIdNoS1} 
 \ee
 where $\widetilde P$ is the  row of the polyphase components of $\widetilde m_0$ without $\widetilde{\mu}_{00}$.
 This modification is useful, since the smaller order of sum rule for $\widetilde m_0$ we demand, 
 the smaller number of non-zero coefficients of $\widetilde m_0$ we obtain. But also,
 the smaller order of smoothness we can get for $\widetilde \phi$.
 And $\widetilde \phi$ may not be in $L_2({\mathbb R}^d).$ So, the desirable order of sum rule for $\widetilde m_0$ should be appropriately chosen.
 \label{remSupp}
    \end{rem}

\subsection{Symmetrization}

The aim of this subsection is to construct the  wavelet masks such that all of them are symmetric in some sense for each matrix from the symmetry group $H$. 
To
do that we extend the definition of $H-$symmetric trigonometric polynomials.
Let $t$ be a trigonometric polynomial. Then $t$ has the $H$-symmetry property if for each matrix $E\in H$ 
 $$t(E^*\xi)= \varepsilon_E
		e^{2\pi i (r_E,\xi)}         
                t(\xi),$$
           where $\varepsilon_E\in {\mathbb C}$, $|\varepsilon_E|=1,$ 
           $r_E\in {\mathbb Z}^{d}.$

Next, we show how to symmetrize  the row of the polyphase components 
of an $H$-symmetric mask, i.e. we want to find a unitary transformation of the row such that each element of the new row has the $H$-symmetry property. This is not always possible. So we restrict our attention to the case when $H$ is an abelian symmetry group on ${\mathbb Z}^{d}$.
	Then  $H$ can be expressed as the direct product of cyclic subgroups due to the fundamental theorem of finite abelian groups.
    Hence, for any number $p=0,\dots,\#\Lambda-1$
    the set $\Gamma_{<s_{p,0}>}$ is an abelian group.
    This is a simple corollary from the fundamental theorem of finite abelian groups since $H_{<s_{p,0}>}$ is a subgroup of  $H$.
    Thus, group $\Gamma_{<s_{p,0}>}$ also can be expressed as the direct product of cyclic subgroups. In order to  write it down, we need some additional notation. Let us fix $p$ and
    let $N_p=\#\Gamma_{<s_{p,0}>}$ and $\gamma_p$ be the number of cyclic subgroups of $\Gamma_{<s_{p,0}>}$.
    Then there exists unique prime numbers $N_{p,i},$ $i=1,\dots,\gamma_p$ and matrices ${\cal E}_1,\dots,{\cal E}_{\gamma_p}\in \Gamma_{<s_{p,0}>}$ such that
    $$\Gamma_{<s_{p,0}>}=\{I_d,{\cal E}_1,\dots,{\cal E}_1^{N_{p,1}-1}\}\times\dots
    \times\{I_d,{\cal E}_{\gamma_p},\dots,{\cal E}_{\gamma_p}^{N_{p,\gamma_p}-1}\}.$$

    Thus, any element $E\in \Gamma_{<s_{p,0}>}$ can be uniquely
     represented as follows $E=\prod_{j=1}^{\gamma_p} {\cal E}_j^{k_j}$, where $k_j$ are some numbers from the sets 
     $\{0,\dots,N_{p,j}-1\},$ respectively.
     Let $L_{p,1}=1,$ $L_{p,i}=\prod_{j=1}^{i-1} N_{p,j},$ for $i=2,\dots,\gamma_p$ and 
     $L_{p,\gamma_p+1}=\prod_{j=1}^{\gamma_p} N_{p,j}=N_p.$
     Thus, each number $k\in\{0,\dots,N_p-1\}$ can be uniquely represented as $k=\sum_{j=1}^{\gamma_p} k_j L_{p,j},$ where $0\le k_j\le N_{p,j}-1.$
     Or equivalently, there is an one-to-one correspondence between $k$ and $(k_1,\dots,k_{\gamma_p}).$ 
     In other words, we get a mixed radix numeral system for all numbers in the set 
     $\{0,\dots,N_{p}-1\}$.
          Let us renumbered all matrices in the group $\Gamma_{<s_{p,0}>}$ according with this numeral system. 
          Let $E\in \Gamma_{<s_{p,0}>},$ 
then $E=\prod_{j=1}^{\gamma_p} {\cal E}_j^{k_j},$
so the matrix $E$ has number $k$, where  $k=\sum_{j=1}^{\gamma_p} k_j L_{p,j}.$ We denote the matrix $E$ with number $k$  by $E^{(k)}.$ 
     
     Let us define  addition on the set  $\{0,\dots,N_p-1\}$ as follows:
     $k\oplus l=\sum_{j=1}^{\gamma_p} (k_j\oplus l_j) L_{p,j},$
     where $k_j\oplus l_j$ is a summation by module $N_{p,j}$,
     namely 
     $k_j\oplus l_j= (k_j+l_j)\mod N_{p,j},$ where $j\in\{1,...,\gamma_p\}.$
    Not hard to see that, 
     $E^{(k)}E^{(l)}=E^{(k \oplus l)}$, where $E^{(k)}, E^{(l)}\in \Gamma_{<s_{p,0}>}.$

Using such notations and assumptions we can rewrite some useful
 properties. Let us fix $p.$ Then the matrix $K\in H$ can be represented as $K=E^{(n)} F, $ where 
  $E^{(n)}\in \Gamma_{<s_{p,0}>}$,
        $F\in H_{<s_{p,0}>}$.
   In this case the map $j(p,i,K)$ simply means $i\oplus n,$
   since $KE^{(i)}=E^{(i\oplus n)} F, $ 
  Thus,~(\ref{fDigK}) can be written as
 \ba
  Ks_{p,i}=
   Mr^{K}_{p,i}+s_{p,i\oplus n}+Kc-c,
   \label{fDigKrevisit}
  \ea
   where $ r_{p,i}^K=M^{-1} 
      E^{(i\oplus n)} M r_{p,0}^{F}.$ 
Let  $m_0$ be an  $H$-symmetric mask with respect to the origin.
  Then conditions~(\ref{fLemSymPoly}),~(\ref{fPolyE_i}) become
        \be
        \mu_{0,p,i}((M^{-1} E^{(n)} M)^*\xi)=\mu_{0,p,i\oplus n}(\xi), \quad  \forall E^{(n)}\in\Gamma_{\langle s_{p,0}\rangle},
        \label{fLemSymPolyNew}
        \ee
  for $i=0,\dots,N_p-1$       and for matrix  $K\in H$ such that 
        $K=E^{(n)} F,$ $E^{(n)}\in \Gamma_{<s_{p,0}>}$,
        $F\in H_{<s_{p,0}>}$
        \be
        \mu_{0,p,i}((M^{-1} K M)^*\xi)=e^{2\pi i (r^K_{p,i},\xi)} \mu_{0,p,i\oplus n}(\xi), 
        \label{fLemSymPolyNew1}
        \ee
 where $ r_{p,i}^K=M^{-1} 
      E^{(i\oplus n)} M r_{p,0}^{F},$ $r_{p,0}^{F}$ is from~(\ref{fDigF}). 

    Let $\varepsilon_{N_{p,i}}= e^{ \frac {2\pi i}{N_{p,i}}}.$    
    For any $p\in\{1,\dots,\#\Lambda-1\}$,
    $W_{N_{p,i}}=\frac 1 {\sqrt{N_{p,i}}}\{
    \varepsilon_{N_{p,i}}^{kl}\}_{k,l=0,N_{p,i}-1}$ is the matrix of the discrete Fourier transform. It is known that $W_{N_{p,i}}$ is a unitary and symmetric matrix, i.e. $W_{N_{p,i}}W_{N_{p,i}}^*=I_m, W_{N_{p,i}}^T=W_{N_{p,i}}.$
Define ${\cal W}_p=W_{N_{p,1}}\otimes\dots\otimes W_{N_{p,\gamma_p}},$ where operation $\otimes$ is the Kronecker product.
This is also a unitary matrix and
the elements can be expressed as follows
$$\left[{\cal W}_p\right]_{k,l}=
\varepsilon_{N_{p,1}}^{k_1 l_1}\dots \varepsilon^{k_{\gamma_p} l_{\gamma_p}}_{N_{p,\gamma_p}},$$
where $k,l$ are corresponded to $(k_1,\dots,k_{\gamma_p})$ and 
$(l_1,\dots,l_{\gamma_p})$ respectively according to the above mixed radix numeral system.

Some of the properties of the matrix ${\cal W}_p$ are $$\left[{\cal W}_p\right]_{k,l}
\left[{\cal W}_p\right]_{n,l}=\left[{\cal W}_p\right]_{k\oplus n,l}, \quad \left[{\cal W}_p\right]_{k,l}\overline{\left[{\cal W}_p\right]_{k,l}}=1, \quad 
k,l,n=0,\dots,N_{p,i}-1.$$

The next Lemma shows that ${\cal W}_p$ symmetrizes the part of the row of polyphase components.

\begin{lem} Let $m_0$ be an $H$-symmetric with respect to the origin mask and $\mu_{0,p,i}$ be its polyphase components, $T_p$ be the row of polyphase components defined by 
$T_p=(\mu_{0,p,0},...,\mu_{0,p,N_p-1}).$
Suppose that $r_{p,0}^F=M^{-1}E Mr_{p,0}^F$ for all  $F\in H_{<s_{p,0}>}$ and $E\in\Gamma_{<s_{p,0}>}.$ Then each element of the row $T_p':=T_p{\cal W}_p$ has the $H$-symmetry property.
\label{propSym}
\end{lem}

{\bf Proof.} Denote the elements of the new row $T_p'$ by $\mu'_{0,p,r}$, $r=0,\dots,N_p-1.$ By definition we have
$\mu'_{0,p,r}(\xi)=
\sum\limits_{k=0}^{N_p-1}\left[{\cal W}_p\right]_{k,r} \mu_{0,p,k}(\xi).$
Firstly, we show that $\mu'_{0,p,r}(\xi)$ has the $(M^{-1}\Gamma_{<s_{p,0}>}M)$-symmetry property.
Indeed, by~(\ref{fLemSymPolyNew}) and the properties of  ${\cal W}_p$  we have for all $E^{(n)}\in\Gamma_{<s_{p,0}>}$
$$\mu'_{0,p,r}((M^{-1}E^{(n)} M)^*\xi)=
\sum\limits_{k=0}^{N_p-1}\left[{\cal W}_p\right]_{k,r}
\mu_{0,p,k}((M^{-1}E^{(n)} M)^*\xi)=
$$
$$
\sum\limits_{k=0}^{N_p-1}\left[{\cal W}_p\right]_{n\oplus k,r}
\overline{\left[{\cal W}_p\right]_{n,r}}
\mu_{0,p,n\oplus k}(\xi)=
\overline{\left[{\cal W}_p\right]_{n,r}} \mu'_{0,p,r}(\xi).
$$
Now show that $\mu'_{0,p,r}(\xi)$, $r=0,\dots,N_p-1,$  are $(M^{-1}H_{<s_{p,0}>}M)$-symmetric.
Indeed, 
	$$\mu'_{0,p,r}((M^{-1}F M)^*\xi)=
	\sum\limits_{k=0}^{N_p-1}\left[{\cal W}_p\right]_{k,r}
\mu_{0,p,0}((M^{-1}F E^{(k)} M)^*(\xi)=
$$$$
\sum\limits_{k=0}^{N_p-1}\left[{\cal W}_p\right]_{k,r}
\mu_{0,p,0}((M^{-1}E^{(k)} M)^*\xi) e^{-2\pi i (M^{-1}E^{(k)}M r_{p,0}^F,\xi)}=
\mu'_{0,p,r}(\xi)e^{-2\pi i (r_{p,0}^F,\xi)}.$$
If $K\in H$ then  $K$ can be represented as $K=E^{(n)} F$, where $F\in H_{<s_{p,0}>}$ and $E^{(n)}\in\Gamma_{<s_{p,0}>}.$  Thus, $\mu'_{0,p,r}(\xi)$ has the $H$-symmetry property
$$\mu'_{0,p,r}((M^{-1}K M)^*\xi)=
\overline{\left[{\cal W}_p\right]_{n,r}}\mu'_{0,p,r}(\xi)e^{-2\pi i (r_{p,0}^{F},\xi)}. \Diamond$$

According with this Lemma we need a special assumption   $r_{p,0}^F=M^{-1}E Mr_{p,0}^F$ for all  $F\in H_{<s_{p,0}>}$ and $E\in\Gamma_{<s_{p,0}>}$   to ensure the $H$-symmetry property for all components of the row $T_p'$.
Notice that this assumption is used only to provide $(M^{-1}H_{<s_{p,0}>}M)$-symmetry for $\mu'_{0,p,r}(\xi)$.

Define a block diagonal unitary matrix ${\cal W}$ as follows:
${\cal W}=diag ({\cal W}_{1},\dots,{\cal W}_{\#\Lambda-1}).$ Let $m_0$ be an $H$-symmetric mask
    with respect to the origin.
       Suppose the special assumption is valid for all 
        $p=1,\dots,\#\Lambda-1,$ i.e.
         $r_{p,0}^F=M^{-1}E Mr_{p,0}^F$ for all  $F\in H_{<s_{p,0}>}$ and $E\in\Gamma_{<s_{p,0}>}.$
    Then the matrix ${\cal W}$ symmetrizes the row $P=(T_1,\dots,T_{
    \#\Lambda-1})$ of the polyphase components of $m_0,$ namely all elements of the row 
$$P^s:=P{\cal W}=(T_1{\cal W}_1,\dots,T_{\#\Lambda-1}{\cal W}_{\#\Lambda-1})$$
    have the $H$-symmetry property. 
  Wavelets with the $H$-symmetry property can be constructed using matrix ${\cal W}$ .   

 \begin{theo} 
  Let $H$ be an abelian symmetry group on ${\mathbb Z}^{d}$, 
  $M$ be an appropriate matrix dilation, $n\in{\mathbb N}$
and let $m_0$ and $\widetilde m_0$ be $H$-symmetric
    with respect to the origin masks such that
     $m_0$
    has sum rule of order $n$,
      $\widetilde m_0$ is constructed by~(\ref{fWM0Def}).
    The special assumption is valid for all 
        $p=1,\dots,\#\Lambda-1,$ i.e.
         $r_{p,0}^F=M^{-1}E Mr_{p,0}^F$ for all  $F\in H_{<s_{p,0}>}$ and $E\in\Gamma_{<s_{p,0}>}.$
    Then there exist wavelet masks  
 which have the $VM^n$ property and the $H$-symmetry property. 
If $\phi,$ $\widetilde\phi \in L_2({\mathbb R}^d)$, then the corresponding 
     wavelet system is a dual wavelet frame.
\label{theoWaveSym}
\end{theo}

{\bf Proof.}  
Let $T_p$ denote the row of the polyphase components defined by 
$T_p=(\mu_{0,p,0},...,
\mu_{0,p,N_p-1}),$
 where $N_p=\#\Gamma_{\langle s_{p,0}\rangle},$ $p=1,...,\#\Lambda-1.$
  Let $P=(T_1,\dots,T_{\#\Lambda-1})$.
Let us consider the matrix extension~(\ref{fNNblockDual}) with $U=\widetilde U={\cal W}^*$
\be
   {\cal N}=\left( \begin{array}{c|c}
      \frac 1 {\sqrt m} & P  \\
      \hline
      &  \\
     - {\cal W}^* P^{*} & \sqrt m( {\cal W}^*- {\cal W}^*P^{*}  P)
    \end{array}\right),
\quad 
{\cal \widetilde N}=\left(
\begin{array}{c|c}
     \sqrt m (1-PP^{*})& P \\
      \hline
      &   \\
      -{\cal W}^* P^{*}  &  \frac 1 {\sqrt m}  {\cal W}^*
    \end{array} \right)
    \ee
    Let us consider submatrix 
    $\sqrt m({\cal W}^{*}-{\cal W}^{*} P^{*}P)$.
    It can be represented as block matrix

\begin{footnotesize}
    \begin{eqnarray}
     {\cal W}^{*}-{\cal W}^{*} P^{*}P=
    \hspace{10.5cm }
    \nonumber
    \\
    \left( \begin{array}{cccc}
      {\cal W}_1^*-(T_1 {\cal W}_1)^* T_1 & 
      	-(T_1 {\cal W}_1)^*  T_2 & \dots & 
      		-(T_1 {\cal W}_1)^*  T_{\#\Lambda-1}\\
      -(T_2 {\cal W}_2)^*  T_1 &
      	{\cal W}_2^*-(T_2 {\cal W}_2)^*  T_2 & \dots & 
      		-(T_2 {\cal W}_2)^*  T_{\#\Lambda-1} \\
      		\vdots & & \ddots & \vdots \\
      	-(T_{\#\Lambda-1} {\cal W}_{\#\Lambda-1})^*  T_1 &	-(T_{\#\Lambda-1} {\cal W}_{\#\Lambda-1})^*  T_2
      	 &\dots &
      	{\cal W}_{\#\Lambda-1}^*-(T_{\#\Lambda-1} 
      	{\cal W}_{\#\Lambda-1})^*  T_{\#\Lambda-1}
      	\nonumber
    \end{array}\right).
    \end{eqnarray}   
\end{footnotesize}
    Denote the elements of the submatrix $\sqrt m({\cal W}^{*}-{\cal W}^{*} P^{*}P)$ as 
    $$\sqrt m({\cal W}^{*}-{\cal W}^{*} P^{*}P):=\{\mu_{(p,i),(t,j)}\}_
    {p=1,\dots,\#\Lambda-1, i=0,\dots,N_p-1}^
    {t=1,\dots,\#\Lambda-1, j=0,\dots,N_t-1}.$$
	With this numeration the element $\mu_{(p,i),(t,j)}$ in the submatrix is the element in the block $(p,t)$ with position $(i,j)$. Namely, if $p\neq t$ then
	$$\mu_{(p,i),(t,j)}(\xi)=\sqrt m\, [-(T_p {\cal W}_p)^*  T_t]_{i,j}=
	-\sqrt m\, \overline{\mu'_{0,p,i}(\xi)} \mu_{0,t,j}(\xi);$$
	if $p=t$ then
	$$\mu_{(p,i),(p,j)}(\xi)=\sqrt m\,[ {\cal W}^*_p-(T_p {\cal W}_p)^*  T_p]_{i,j}=\sqrt m\,\left(
	\overline{[ {\cal W}_p]_{i,j}}-\overline{\mu'_{0,p,i}(\xi)} \mu_{0,p,j}(\xi)\right),$$
	where $\mu'_{0,p,i},$ $ i=0,\dots,N_p-1$, are the elements of the row $T_p {\cal W}_p$,
   $ p=1,\dots,\#\Lambda-1,$. All $\mu'_{0,p,i}$ have the $H$-symmetry property by Lemma~\ref{propSym}.
	Now, we check the $H-$symmetry property of $\mu_{(p,i),(t,j)}$. Let us fix $p$ and $t$, $p\neq t$ and  matrix $K\in H.$ Then we can represent  $K$ as $K=E_p^{(k)}F_p$, where $E^{(k)}_p\in \Gamma_{<s_{p,0}>}$,
        $F_p\in H_{<s_{p,0}>},$ 
        and as $K=E_t^{(l)}F_t,$ where
        $E^{(l)}_t\in \Gamma_{<s_{t,0}>}$,
        $F_t\in H_{<s_{t,0}>}.$  Therefore, by ~(\ref{fLemSymPolyNew}),~(\ref{fLemSymPolyNew1}) and Lemma~\ref{propSym} and we have
        
	\be
	\mu_{(p,i),(t,j)}((M^{-1}K M)^*\xi)=
	-\sqrt m\, \overline{\mu'_{0,p,i}((M^{-1}K M)^*\xi)}
	\mu_{0,t,j}((M^{-1}K M)^*\xi)=$$
	$$-\sqrt m\, [ {\cal W}_p]_{i,k}
	e^{2\pi i (r^{F_p}_{p,0},\xi)} 
	\overline{\mu'_{0,p,i}(\xi)}
	e^{-2\pi i (r^{F_t}_{t,0},\xi)}
	 \mu_{0,t,j\oplus l}(\xi) 
	=$$
	$$[ {\cal W}_p]_{i,k} 
	e^{2\pi i (r^{F_p}_{p,0},\xi)} 
		e^{-2\pi i (r^{F_t}_{t,0},\xi)}\mu_{(p,i),(t,j\oplus l)}(\xi)
	.	
	\label{fTheoHSym}
	\ee
		Let $p=t$ and $K\in H.$ Then matrix $K$ can be represented as $K=E_p^{(k)}F_p$, where $E^{(k)}_p\in \Gamma_{<s_{p,0}>}$,
        $F_p\in H_{<s_{p,0}>}$. Therefore,
$$
        \mu_{(p,i),(p,j)}((M^{-1}K M)^*\xi)=
        \sqrt m\,\left(
	\overline{[ {\cal W}_p]_{i,j}}-
	\overline{\mu'_{0,p,i}((M^{-1}K M)^*\xi)}
	\mu_{0,p,j}((M^{-1}K M)^*\xi)\right)=$$
	$$\sqrt m\,\left(\overline{[ {\cal W}_p]_{i,j}}-
	[ {\cal W}_p]_{i,k}
	e^{2\pi i (r^{F_p}_{p,0},\xi)} 
	\overline{\mu'_{0,p,i}(\xi)}
	e^{-2\pi i (r^{F_p}_{p,0},\xi)}
	\mu_{0,p,j\oplus k}(\xi)\right)=$$
	$$\sqrt m\,[ {\cal W}_p]_{i,k}(\overline{[ {\cal W}_p]_{i,j\oplus k}}-
	\overline{\mu'_{0,p,i}(\xi)}
	\mu_{0,p,j\oplus k}(\xi))=
	[ {\cal W}_p]_{i,k}\mu_{(p,i),(p,j\oplus k)}(\xi).
$$
   The elements of the first column $-{\cal W}^* P^{*}$ in matrix ${\cal N}$ we denote by 
   $\mu_{(p,i),{(0,0)}}$ and they are equal to 
   $\mu_{(p,i),{(0,0)}}=-\frac{1}{\sqrt{m}}\overline{\mu_{0,p,i}'(\xi)}$. Since 
   $\Gamma_{<s_{0,0}>}=\{I_d\}$ and $r_{0,0}^F={\bf0},$ $\forall F\in H_{<s_{0,0}>}$, equation~(\ref{fTheoHSym}) for $\mu_{(p,i),{(0,0)}}$
   is as follows 
   $$
   \mu_{(p,i),(0,0)}((M^{-1}K M)^*\xi)=
	[ {\cal W}_p]_{i,k} 
	e^{2\pi i (r^{F_p}_{p,0},\xi)} 
		\mu_{(p,i),(0,0)}(\xi),
		\quad\forall  K\in H.
   $$
   Let us fix numbers $p$ and $i$ and construct the wavelet mask $m_{(p,i)}$ from the trigonometric polynomials $\mu_{(p,i),(t,j)}$  by~(\ref{PR}).
Let us show that the resulting wavelet  mask $m_{(p,i)}$ has the $H$-symmetry
 property. Suppose $K\in H.$  Then we can represent matrix $K$ as $K=E_p^{(k)}F_p$, where $E^{(k)}_p\in \Gamma_{<s_{p,0}>}$,
        $F_p\in H_{<s_{p,0}>},$ and $K=E_t^{(l)}F_t,$ where
        $E^{(l)}_t\in \Gamma_{<s_{t,0}>}$,
        $F_t\in H_{<s_{t,0}>}.$ 
 Using~(\ref{fDigKrevisit}) and~(\ref{fTheoHSym}), we get
 $$m_{(p,i)}(K^*\xi)= \frac1{\sqrt m}\sum\limits_{t=0}^{\#\Lambda-1} \sum\limits_{j=0}^{
 \#\Gamma_{<s_{t,0}>}-1} 
 e^{2\pi i(K s_{t,j},\xi)}\mu_{(p,i),(t,j)}(M^*K^*\xi)=$$
                $$\frac { [ {\cal W}_p]_{i,k} }{\sqrt m}               
                \sum\limits_{t=0}^{\#\Lambda-1} 
                \sum\limits_{j=0}^{
 \#\Gamma_{<s_{t,0}>}-1} 
                 e^{2\pi i(s_{t,j\oplus l},\xi)}
                 e^{2\pi i (Mr_{t,j}^K,\xi)}            
                \mu_{(p,i),(t,j\oplus l)}(M^*\xi)
	e^{2\pi i (Mr^{F_p}_{p,0}-Mr^{F_t}_{t,0},\xi)} 
	=$$
		$$       [ {\cal W}_p]_{i,k}		
                e^{2\pi i (Mr^{F_p}_{p,0},\xi)} 
                m_{(p,i)}(\xi).$$
                The last equality is valid since  $r_{t,j}^K=M^{-1}E_t^{(j\oplus l)}M r_{t,0}^{F_t}=r_{t,0}^{F_t}.$ 
                Thus, we get the $H$-symmetry property for 
                 the wavelet mask $m_{(p,i)}$.
         In our case the $H$-symmetry
 property means that when $K=F_p,$ then  $ m_{(p,i)}(\xi)$ is 
        $H_{<s_{p,0}>}$-symmetric with respect to the center $s_{p,0}.$
        When $K=E^{(k)}_p$ then $ m_{(p,i)}(\xi)$ has the $\Gamma_{<s_{p,0}>}$-symmetry property, namely,
        $m_{(p,i)}(E^{(k)*}_p\xi)=[ {\cal W}_p]_{i,k}		\,\,         
                m_{(p,i)}(\xi).$
                Not hard to see that the dual wavelet masks $\widetilde m_{(p,i)}$ have the same $H$-symmetry property as $m_{(p,i)}$.
                
                     The $VM^n$ property for the  wavelet masks is provided by Lemma~\ref{lemSkVM}.
                     $\Diamond$
                     
Note that the comments in Remark~\ref{remSupp} are also true for Theorem~\ref{theoWaveSym}.


\section{Examples}

In this section we
     give several examples to illustrate the main results of this paper.
    All examples are based on the construction of 
   $H$-symmetric interpolatory refinable masks $m_0$ and $\widetilde m_0$ and wavelet masks $m_{\nu}, $ $\widetilde m_{\nu},$ $\nu=1,\dots,m$ by Theorems~\ref{theoMaskH},~\ref{theoWaveNoSym} and~\ref{theoWaveSym}.

    1. Let $H$ be a hexagonal symmetry group, namely,
  \begin{small}
    
        $$H=\left\{\pm I_2,\pm\left(
\begin{array}{cc}
 0 & 1 \\
 1 & 0 \\
\end{array}
\right),\pm\left(
\begin{array}{cc}
 1 & 0 \\
 1 & -1 \\
\end{array}
\right),\pm\left(
\begin{array}{cc}
 1 & -1 \\
 1 & 0 \\
\end{array}
\right),\pm\left(
\begin{array}{cc}
 0 & 1 \\
 -1 & 1 \\
\end{array}
\right),\pm\left(
\begin{array}{cc}
 -1 & 1 \\
 0 & 1 \\
\end{array}
\right)\right\}$$
  \end{small}
         $c={\bf0},$ ${M}= \left(\begin{matrix} 2&-1\cr
    1& 1 \cr
    \end{matrix}\right).
    $
        The set of digits is $D(M)=\{s_0=(0,0), s_1=(0,1), s_2=(0,-1)\},$ $m=3.$ 
        Let us construct a refinable mask  that is $H$-symmetric with respect to the origin
        and has sum rule of order $n=3.$
        In our case, the digits are renumbered as $s_{0,0}=(0,0),$ $s_{1,0}=(0,1),s_{1,1}=(0,-1).$
       And $H_{<s_{0,0}>}=H,$ $\Gamma_{<s_{0,0}>}=\{I_2\},$
       $$H_{<s_{1,0}>}=\left\{I_2,\left(
\begin{array}{cc}
 0 & 1 \\
 1 & 0 \\
\end{array}
\right),\left(
\begin{array}{cc}
 -1 & 0 \\
 -1 & 1 \\
\end{array}
\right),\left(
\begin{array}{cc}
 -1 & 1 \\
 -1 & 0 \\
\end{array}
\right),\left(
\begin{array}{cc}
 0 & -1 \\
 1 & -1 \\
\end{array}
\right),\left(
\begin{array}{cc}
 1 & -1 \\
 0 & -1 \\
\end{array}
\right)\right\},$$
$\Gamma_{<s_{0,0}>}=\left\{I_2,\left(
\begin{array}{cc}
 1 & 0 \\
 1 & -1 \\
\end{array}
\right)\right\}.$
        According to Theorem~\ref{theoMaskH},
        the mask $m_0$
        can be constructed as follows

            $$m_0: \left(
\begin{array}{ccccc}
 0 & 0 & -\frac{1}{27} & 0 & -\frac{1}{27} \\
 0 & 0 & \frac{4}{27} & \frac{4}{27} & 0 \\
 -\frac{1}{27} & \frac{4}{27} & \frac{\textbf{1}}{\textbf{3}} & \frac{4}{27} & -\frac{1}{27} \\
 0 & \frac{4}{27} & \frac{4}{27} & 0 & 0 \\
 -\frac{1}{27} & 0 & -\frac{1}{27} & 0 & 0 \\
\end{array}
\right)$$
    with support in $[-2,2]^2\bigcap{\mathbb Z}^2.$ Note that $\mu_{00}=\frac 1 {\sqrt 3}.$
        The corresponding refinable function $\phi$ is in $L_2({\mathbb R}^2)$ and $\nu_2(\phi)\ge   1.8959.$
        The Sobolev smothness exponent is calculated by~\cite[Theorem 7.1]{HanSymSmooth}.
       The dual mask $\widetilde m_0 $ is
               $$\widetilde m_0: \left(
\begin{array}{ccccccccc}
 0 & 0 & 0 & 0 & 0 & 0 & -\frac{2}{243} & 0 & 0 \\
 0 & 0 & 0 & 0 & \frac{8}{243} & 0 & 0 & \frac{8}{243} & 0 \\
 0 & 0 & -\frac{2}{243} & 0 & -\frac{1}{27} & -\frac{16}{243} & -\frac{1}{27} & 0 & -\frac{2}{243}
   \\
 0 & 0 & 0 & -\frac{16}{243} & \frac{4}{27} & \frac{4}{27} & -\frac{16}{243} & 0 & 0 \\
 0 & \frac{8}{243} & -\frac{1}{27} & \frac{4}{27} & \frac{\textbf{47}}{\textbf{81}} & \frac{4}{27} & -\frac{1}{27} &
   \frac{8}{243} & 0 \\
 0 & 0 & -\frac{16}{243} & \frac{4}{27} & \frac{4}{27} & -\frac{16}{243} & 0 & 0 & 0 \\
 -\frac{2}{243} & 0 & -\frac{1}{27} & -\frac{16}{243} & -\frac{1}{27} & 0 & -\frac{2}{243} & 0 & 0
   \\
 0 & \frac{8}{243} & 0 & 0 & \frac{8}{243} & 0 & 0 & 0 & 0 \\
 0 & 0 & -\frac{2}{243} & 0 & 0 & 0 & 0 & 0 & 0 \\
\end{array}
\right)$$   
 with support in $[-4,4]^2\bigcap{\mathbb Z}^2.$
        The corresponding refinable function $\widetilde\phi$ is in $L_2({\mathbb R}^2)$ and $\nu_2(\widetilde\phi)\ge   0.7852.$ 
 Since the  hexagonal symmetry group  $H$ is not abelian, the wavelet masks are constructed by Theorem~\ref{theoWaveNoSym}:

            $$m_1: \left(
\begin{array}{ccccccccc}
 0 & 0 & 0 & 0 & -\frac{1}{243} & 0 & -\frac{1}{243} & 0 & 0 \\
 0 & 0 & 0 & 0 & \frac{4}{243} & \frac{8}{243} & 0 & \frac{4}{243} & 0 \\
 0 & 0 & -\frac{1}{243} & \frac{8}{243} & \frac{1}{27} & -\frac{8}{243} & -\frac{2}{27} & 0 &
   -\frac{1}{243} \\
 0 & 0 & 0 & -\frac{8}{243} & -\frac{8}{81} & -\frac{4}{27} & -\frac{8}{243} & \frac{8}{243} & 0 \\
 0 & \frac{4}{243} & -\frac{2}{27} & -\frac{4}{27} & \frac{64}{81} & -\frac{8}{81} & \frac{1}{27} &
   \frac{4}{243} & -\frac{1}{243} \\
 0 & 0 & -\frac{8}{243} & -\frac{8}{81} & -\frac{\textbf{4}}{\textbf{27}} & -\frac{8}{243} & \frac{8}{243} & 0 & 0 \\
 -\frac{1}{243} & \frac{8}{243} & \frac{1}{27} & -\frac{8}{243} & -\frac{2}{27} & 0 &
   -\frac{1}{243} & 0 & 0 \\
 0 & \frac{4}{243} & \frac{8}{243} & 0 & \frac{4}{243} & 0 & 0 & 0 & 0 \\
 -\frac{1}{243} & 0 & -\frac{1}{243} & 0 & 0 & 0 & 0 & 0 & 0 \\
\end{array}
\right),$$

$$\widetilde m_1:
\left(
\begin{array}{ccccc}
 0 & 0 & \frac{1}{27} & 0 & 0 \\
 0 & 0 & 0 & -\frac{4}{27} & 0 \\
 0 & -\frac{4}{27} & \frac{1}{3} & 0 & \frac{1}{27} \\
 0 & 0 & -\frac{\textbf{4}}{\textbf{27}} & 0 & 0 \\
 \frac{1}{27} & 0 & 0 & 0 & 0 \\
\end{array}
\right), \quad m_2=m_1(E^*\xi), \quad \widetilde m_2=\widetilde m_1(E^*\xi),$$
where $E=\left(
\begin{array}{cc}
 1 & 0 \\
 1 & -1 \\
\end{array}
\right).$ The bold element in the matrices corresponds to the coefficient $h_{{\bf0}}$ of the masks.
%
 The wavelet masks $ m_1, m_2$, $\widetilde m_1, \widetilde m_2$ have vanishing moments of order 3, the masks $ m_1, \widetilde m_1$ are $H_{<s_{1,0}>}$-symmetric. 

    2. Let $H=H^{axis}$,
         $M= \left(\begin{matrix} 3&0\cr
    0& 2 \cr
    \end{matrix}\right).
    $
        Let us construct a refinable mask $m_0$ that is $H$-symmetric with respect to the origin
        and has sum rule of order $n=1.$
         In our case, the digits are renumbered as  $s_{0,0}=(0,0),$ $s_{1,0}=(0,1)$, 
        $s_{2,0}=(1,0), s_{2,1}=(-1,0)$,
            $s_{3,0}=(1,1), s_{3,1}=(-1,1)$,
       $H_{<s_{0,0}>}=H_{<s_{1,0}>}=H,$ $\Gamma_{<s_{0,0}>}=\Gamma_{<s_{1,0}>}=\{I_2\},$
       $$H_{<s_{2,0}>}=H_{<s_{3,0}>}=
       \left\{I_2,\left(
\begin{array}{cc}
 1 & 0 \\
 0 & -1 \\
\end{array}
\right)\right\},\quad \Gamma_{<s_{2,0}>}=
\Gamma_{<s_{3,0}>}=\left\{I_2,\left(
\begin{array}{cc}
 -1 & 0 \\
 0 & 1 \\
\end{array}
\right)\right\}.$$

       The refinable mask $m_0$ 
constructed by Theorem~\ref{theoMaskH} and  the dual mask $\widetilde m_0$  constructed by~(\ref{fWM0Def}) are

            $$m_0: 
           \left(
\begin{array}{ccc}
 \frac{1}{12} & \frac{1}{12} & \frac{1}{12} \\
 \frac{1}{6} & \frac{\textbf{1}}{\textbf{6}} & \frac{1}{6} \\
 \frac{1}{12} & \frac{1}{12} & \frac{1}{12} \\
\end{array}
\right), \quad 
\widetilde m_0: 
               \left(
\begin{array}{ccc}
 0 & -\frac{1}{8} & 0 \\
 \frac{1}{12} & \frac{1}{12} & \frac{1}{12} \\
 \frac{1}{6} & \frac{\textbf{5}}{\textbf{12}} & \frac{1}{6} \\
 \frac{1}{12} & \frac{1}{12} & \frac{1}{12} \\
 0 & -\frac{1}{8} & 0 \\
\end{array}
\right),$$   
    with support in $[-1,1]^2\bigcap{\mathbb Z}^2$ and in $[-1,1]\times[-2,2]\bigcap{\mathbb Z}^2$ respectively.
 Since the  axial symmetry group is an abelian group, the wavelet masks can be constructed by Theorem~\ref{theoWaveSym}. Note that the special assumption is valid in this case. To avoid roots in  the coefficients of wavelet mask   we take matrix ${\cal W}$ as
 $${\cal W}=
 \left(
\begin{array}{ccccc}
 1 & 0 & 0 & 0 & 0 \\
 0 & 1 & 1 & 0 & 0 \\
 0 & 1 & -1 & 0 & 0 \\
 0 & 0 & 0 & 1 & 1 \\
 0 & 0 & 0 & 1 & -1 \\
\end{array}
\right), \quad 
\widetilde {\cal W}=
\left(
\begin{array}{ccccc}
 1 & 0 & 0 & 0 & 0 \\
 0 & \frac{1}{2} & \frac{1}{2} & 0 & 0 \\
 0 & \frac{1}{2} & -\frac{1}{2} & 0 & 0 \\
 0 & 0 & 0 & \frac{1}{2} & \frac{1}{2} \\
 0 & 0 & 0 & \frac{1}{2} & -\frac{1}{2} \\
\end{array}
\right),$$
where $\widetilde {\cal W}$ is a paraunitary matrix for ${\cal W}$,
i.e. ${\cal W}^* \widetilde {\cal W}
=I_5$. 
The wavelet masks $m_1, m_2, m_3, m_4, m_5$ are
            $$
          \begin{footnotesize}
          \left(
\begin{array}{ccc}
 -\frac{1}{24} & -\frac{1}{24} & -\frac{1}{24} \\
 -\frac{1}{12} & -\frac{1}{12} & -\frac{1}{12} \\
 -\frac{1}{12} & \frac{11}{12} & -\frac{1}{12} \\
 -\frac{1}{12} & -\frac{\textbf{1}}{\textbf{12}} & -\frac{1}{12} \\
 -\frac{1}{24} & -\frac{1}{24} & -\frac{1}{24} \\
\end{array}
\right),\left(
\begin{array}{ccc}
 -\frac{1}{6} & -\frac{1}{6} & -\frac{1}{6} \\
 \frac{2}{3} & -\frac{\textbf{1}}{\textbf{3}} & \frac{2}{3} \\
 -\frac{1}{6} & -\frac{1}{6} & -\frac{1}{6} \\
\end{array}
\right),\left(
\begin{array}{ccc}
 0 & 0 & 0 \\
 -1 & \textbf{0} & 1 \\
\end{array}
\right),
          \end{footnotesize}$$
$$\begin{footnotesize}
\left(
\begin{array}{ccc}
 -\frac{1}{12} & -\frac{1}{12} & -\frac{1}{12} \\
 -\frac{1}{6} & -\frac{1}{6} & -\frac{1}{6} \\
 \frac{5}{6} & -\frac{1}{6} & \frac{5}{6} \\
 -\frac{1}{6} & -\frac{\textbf{1}}{\textbf{6}} & -\frac{1}{6} \\
 -\frac{1}{12} & -\frac{1}{12} & -\frac{1}{12} \\
\end{array}
\right),\left(
\begin{array}{ccc}
 -1 & 0 & 1 \\
 0 & \textbf{0} & 0 \\
\end{array}
\right)
\end{footnotesize}
            ,$$
 The dual wavelet masks $\widetilde m_1,\widetilde m_2,\widetilde m_3,\widetilde m_4,\widetilde m_5$ are          
            $$
           \begin{footnotesize}
            \left(
\begin{array}{ccc}
 0 & -\frac{1}{12} & 0 \\
 0 & \frac{1}{6} & 0 \\
 0 & -\frac{\textbf{1}}{\textbf{12}} & 0 \\
\end{array}
\right),\left(
\begin{array}{ccc}
 0 & 0 & 0 \\
 \frac{1}{12} & -\frac{\textbf{1}}{\textbf{6}} & \frac{1}{12} \\
\end{array}
\right),\left(
\begin{array}{ccc}
 0 & 0 & 0 \\
 -\frac{1}{12} & \textbf{0} & \frac{1}{12} \\
\end{array}
\right),
           \end{footnotesize}
           $$
$$
\begin{footnotesize}
\left(
\begin{array}{ccc}
 0 & -\frac{1}{12} & 0 \\
 \frac{1}{12} & 0 & \frac{1}{12} \\
 0 & -\frac{\textbf{1}}{\textbf{12}} & 0 \\
\end{array}
\right),\left(
\begin{array}{ccc}
 -\frac{1}{12} & 0 & \frac{1}{12} \\
 0 & \textbf{0} & 0 \\
\end{array}
\right),
\end{footnotesize}
            $$
where the bold element in the matrices corresponds to the coefficient $h_{{\bf0}}$ of the  masks.
All wavelet masks  have vanishing moments of order 1 and the $H^{axis}-$symmetry property.

    \end{document}